\documentclass{article}[12pt]

\usepackage{amssymb, amsmath, amsthm,commath,stackrel,enumerate,mathtools,bbold,xcolor,thm-restate,comment, geometry}
\usepackage{stmaryrd}
\usepackage{amsfonts}
\usepackage{amssymb}
\usepackage{amsthm}
\usepackage{amsmath}
\usepackage{amscd}
\usepackage{hyperref}
\usepackage{tikz}
\usepackage{pgf}
\usepackage{graphicx}

\renewcommand\labelenumi{(\roman{enumi})}
\renewcommand\theenumi\labelenumi
\DeclareMathOperator{\OB}{\mathcal O\mathcal B}
\DeclareMathOperator{\cont}{cont}

\author{Romain Tessera, Jeroen Winkel}
\date{}
\title{Coarse fixed point properties}

\newtheorem{thmIntro}{Theorem}[]

\theoremstyle{definition}
\newtheorem{example}{Example}[section]
\newtheorem{examples}{Examples}[section]
\theoremstyle{plain}
\newtheorem{theorem}[example]{Theorem}
\newtheorem{lemma}[example]{Lemma}
\newtheorem{construction}[example]{Construction}
\newtheorem{proposition}[example]{Proposition}

\newtheorem{corollary}[example]{Corollary}

\newtheorem{que}[example]{Question}
\theoremstyle{remark}
\newtheorem{remark}[example]{Remark}
\theoremstyle{definition}
\newtheorem{definition}[example]{Definition}

\renewcommand{\phi}{\varphi}
\renewcommand{\H}{\mathcal H}
\DeclareMathOperator{\N}{\mathbb N}
\DeclareMathOperator{\R}{\mathbb R}
\DeclareMathOperator{\Z}{\mathbb Z}
\renewcommand{\epsilon}{\varepsilon}
\DeclareMathOperator{\cs}{cs}
\DeclareMathOperator{\diam}{diam}
\DeclareMathOperator{\tensor}{\otimes}
\DeclareMathOperator{\Aut}{Aut}
\DeclareMathOperator{\Homeo}{Homeo}
\DeclareMathOperator{\Isom}{Isom}
\DeclareMathOperator{\Cay}{Cay}
\DeclareMathOperator{\co}{co}
\DeclareMathOperator{\id}{id}
\DeclareMathOperator{\cco}{\overline{co}}

\DeclareMathOperator{\acts}{\curvearrowright}
\DeclareMathOperator{\eps}{\varepsilon}
\DeclareMathOperator{\born}{born}
\begin{document}
\maketitle

\begin{abstract}
We investigate fixed point properties for isometric actions of topological groups on a wide class of metric spaces, with a particular emphasis on Hilbert spaces. Instead of requiring the action to be continuous, we assume that it is ``controlled", i.e. compatible with respect to some natural left-invariant coarse structure. For locally compact groups, we prove that these coarse fixed point properties are equivalent to the usual ones, defined for continuous actions. We deduce generalisations of two results of Gromov originally stated for discrete groups. 
For Polish groups with bounded geometry (in the sense of Rosendal), we prove a version of Serre's theorem on the stability of coarse property FH under central extensions. As an application we prove that the group $\Homeo^+_{\mathbb Z}(\mathbb R)$ has property FH. Finally, we characterise geometric property (T) for sequences of finite Cayley graphs in terms of coarse property FH of a certain group.
\end{abstract}
\section{Introduction}
In this article we investigate geometric properties of {\it bornological} groups, namely groups equipped with a left-invariant coarse structure in the sense of Roe (see Definition \ref{def:bornological}). 
More specifically, we shall focus on {\it fixed point properties} for isometric actions of bornological groups on various classes of metric spaces, with a special emphasis on Hilbert spaces.  
Bornological groups arise naturally in different contexts, notably in the context of topological groups. Indeed, given a topological group, one can consider two natural left-invariant coarse structures: one for which bounded sets are relatively compact subsets, and another one introduced by Rosendal, for which bounded sets are those that are bounded relative to any left-invariant continuous pseudo-metric \cite{Ros21}. For $\sigma$-compact locally compact groups, these two coarse structures coincide.

There is a natural notion of {\it controlled} isometric action of a bornological groups on a metric space  (see Definition \ref{def:CoarseFX}). Given a class $\mathcal X$ of metric spaces, we say that a bornological group $G$ has {\it  coarse} property F$\mathcal X$ if every controlled action of $G$ on a space from $\mathcal X$ admits a fixed point.
For a topological group $G$ (which we always assume to be Hausdorff), we shall say that $G$ has {\it topological} property F$\mathcal X$ if every {\it continuous} isometric action of $G$ on any space in $\mathcal X$ has a fixed point.

We shall work with metric spaces satisfying a certain convexity property: \emph{uniformly convex} metric spaces  (see Definition \ref{def:uniformly convex}). This is motivated by the fact that such spaces admit a natural notion of centre for bounded sets (see Lemma \ref{lem-of-centre}). This class encompasses a wide family of complete geodesic metric spaces such as uniformly convex (real) Banach spaces (in particular $L^p$-spaces for $1<p<\infty$) and CAT(0) spaces. 
We recall that by the Mazur-Ulam theorem, isometric actions on Banach spaces are affine. 
  
This introduction is organised as follows.
\begin{itemize}

\item   First, in \S \ref{secIntro: Borno}, we define bornological groups and discuss some first examples. We also define a coarse version of Kazhdan's property (T), and we state the relevant version of the Delorme-Guichardet theorem relating coarse property (T) and coarse property FH.

\item  In \S \ref{secIntro:LocCompact}, we state our main result regarding locally compact groups, which roughly says that topological and coarse fixed point properties coincide. We use it to extend two theorems, that were originally proved by Gromov for discrete groups, to locally compact groups. 

\item  In \S \ref{secIntro:BoundedGeom}, we consider bornological groups with {\it bounded geometry}. Besides locally compact groups this class contains interesting examples of Polish groups studied by Rosendal \cite{Ros21}. We prove an analogue of Serre's theorem on the stability of property (T) under central extensions for bornological groups with bounded geometry.
We apply this to prove that certain Polish groups which are quasi-isometric to $\Z$, such as $\Homeo^+_{\mathbb Z}(\mathbb R)$, have coarse (and therefore topological) property FH.

\item In \S \ref{secIntro:BoundedProducts}, we consider a topological group naturally associated to a sequence of Cayley graphs of finite groups. We show that {\it topological} and {\it coarse} property (T) of this group encode {\it expansion} and {\it geometric property (T)} (in the sense of Willett and Yu) respectively. 
\end{itemize}

 \subsection{Bornological groups}\label{secIntro: Borno}  

We recall the following notion essentially due to Roe \cite{Roe} and studied by Rosendal in the context of Polish groups \cite{Ros21}.
\begin{definition}\label{def:bornological}
Let $G$ be a group.
A \emph{bornology} is a subset $\mathcal B\subseteq P(G)$ of so-called \emph{bounded sets}, satisfying the following conditions:
\begin{enumerate}
    \item For each $g\in G$, the singleton set $\{g\}$ is bounded.
    \item If $A\subseteq B\subseteq G$ and $B$ is bounded, then $A$ is bounded as well.
    \item If $A\subseteq G$ is bounded, $A^{-1}=\{a^{-1}\mid a\in A\}$ is bounded as well.
    \item If $A,B\subseteq G$ are bounded, then $A\cup B$ and $A\cdot B$ are bounded.
\end{enumerate}
A \emph{bornological group} is a group equipped with a bornology.
\end{definition}

Given a bornological group $(G,\mathcal B)$ we can give $G$ a coarse structure in two different ways: there is the left coarse structure 
\[\mathcal E_L = \left\{E\subseteq G\times G\mid \{{x^{-1}y\mid (x,y)\in E}\}\text{ is bounded}\right\}.\]
and there is the right coarse structure
\[\mathcal E_R = \left\{E\subseteq G\times G\mid \{{xy^{-1}\mid (x,y)\in E}\}\text{ is bounded}\right\}.\]
The left coarse structure is left invariant, and each controlled set is contained in a set of the form $E_A = \{(ag,a)\mid g\in G,a\in A\}$ for some bounded $A\subseteq G$.
A similar formula holds for the right coarse structure.
In general, the left coarse structure and the right coarse structure might not be the same.
In this case, the multiplication is not a controlled map $G\times G\to G$.

\begin{remark}
This notion of bornological group exactly captures the group objects in the category of bornological spaces with controlled maps.
Even though bornological groups also have a coarse structure, we refrain from calling them coarse groups, exactly because they do not correspond to group objects in the category of coarse spaces with controlled maps -- this is because the multipication map does not need to be controlled.
\end{remark}

\begin{examples}Let us list a few examples.
\begin{itemize}
    \item On each group $G$, there are two trivial bornologies: the \emph{discrete bornology} consisting of all finite subsets of $G$, and the bornology containing every subset of $G$.
    \item Let $\sigma: G\acts X$ be an isometric action of a group $G$ on a metric space $X$. Then $\mathcal B_\sigma$ consisting of all subsets $A$ such that $\sigma(A)\cdot x$ has finite diameter for all $x\in X$ is a bornology. 
    \item 
    Let $G$ be a topological group.
Then we can make it a bornological group with the bornology $\mathcal K = \{A\subseteq G\mid A\text{ is precompact}\}$.
We can also make it a bornological group with the bornology $\OB = \cap_\sigma \mathcal B_\sigma$,
where $\sigma$ runs over all continuous isometric actions on metric spaces.
This bornology has been introduced in \cite{Ros21} (see Definition \ref{def:OB} for more details).
For sigma-compact locally compact groups, $\mathcal K = \OB$.

\item Let $(X,\mathcal E)$ be a coarse space.
Let $\Aut_0(X)$ denote the group of bijections $g:X\to X$ for which the graph $\{(x,gx)\mid x\in X\}$ is a controlled set.
We give this a bornology as follows: the set $A\subseteq \Aut_0(X)$ is bounded if the set $\{(x,ax)\mid x\in X,a\in A\}$ is a controlled set.
\end{itemize}
\end{examples}
\begin{definition}\label{def:CoarseFX}
Let $(G,\mathcal B)$ be a bornological group. 
An isometric action $\sigma: G\acts X$ is called controlled if $\mathcal B\subseteq \mathcal B_\sigma$, i.e.\  for $A\in \mathcal B$, $\sigma(A)$ has bounded orbits.
 A bornological group has coarse property F$X$ if all its controlled isometric actions on $X$ have a fixed point.
\end{definition}

Note that a norm-preserving representation on a Banach space is automatically controlled (since orbits are obviously bounded).

We now proceed to define a convenient notion of property (T) for bornological groups. We recall that all unitary representations are controlled, so we will not impose any restriction on them. However, we need to specify the meaning of ``almost'' invariant vectors. For simplicity, we focus on Hilbert spaces in the introduction (see Definition \ref{def:T_B} for a more general treatment).
 
\begin{definition} Let $G$ be a bornological group.
A representation $\pi\colon G\to U(\H)$ has \emph{almost invariant vectors} if for each $\epsilon>0$ and bounded $A\subseteq G$, there is a unit vector $v\in \H$ such that $\sup_{a\in A}\norm{\pi(a)v-v}<\epsilon$.
\end{definition}

\begin{definition}
A bornological group $G$ has \emph{coarse property (T)} if every representation with almost invariant vectors has an invariant vector.
\end{definition}

Adapting the proof of Delorme and Guichardet's famous theorem, we prove that coarse property (T) implies coarse property FH, and the converse is true if the bornology is generated by a single bounded set (see Theorem \ref{prop-Delorme-Guichardet}).

\subsection{The case of locally compact groups}\label{secIntro:LocCompact}
The following result says that under a mild assumption, when studying isometric actions of topological groups on uniformly convex metric spaces, one can restrict to continuous ones. 
\begin{proposition}[see \S \ref{sec:prop:continuousSubaction}]\label{prop:continuousSubaction}
Let $G$ be a topological group and let $U$ be an open neighbourhood of the identity. Let $X$ be a uniformly convex metric space.
Let $\alpha\colon G\curvearrowright X$ be an isometric action.
Let $c>0$ and suppose there is $x\in X$ such that $d(\alpha(g)x,x)\leq c$ for all $g\in U$.
Then there is $y\in X$ with $d(x,y)\leq 2c$, such that the restriction of $\alpha$ to the closed convex hull of $\alpha(G)y$ is continuous.
\end{proposition}
If $X$ is a Banach space, then the closed convex hull of $\alpha(G)y$ can be replaced by the closed vector subspace spanned by $\alpha(G)y$ in the proposition.
We refer to the discussion in \S \ref{sec:questions} for a cohomological  interpretation of Proposition \ref{prop:continuousSubaction} in the case where $X$ is a Banach space.

We give two different proofs of this result, both based on certain notions of centres (see \S \ref{sec:centres}). One of them only applies to actions of locally compact groups on uniformly convex Banach spaces. However, this one exploits a new notion of centre of a finitely additive measure, which we believe could be useful for future applications.

Proposition \ref{prop:continuousSubaction} has the following consequence, formulated in terms of bornological groups.

\begin{thmIntro}[See Corollary \ref{cor:Tdense}]\label{thm:topo/coarseFX}
Let $X$ be a uniformly convex metric space. Equipped with the bornology $\mathcal K$, a locally compact group $G$ has topological property F$X$ if and only if it has coarse property F$X$.
\end{thmIntro}
Note that coarse property F$X$ is obviously stronger than topological property F$X$, hence Theorem \ref{thm:topo/coarseFX} is really about the other implication.
Similar statements can be derived for properties defined in terms of norm-preserving representations (see Corollary \ref{cor:THetc}).

We now present two consequences, which are new in general for locally compact groups, but are due to Gromov \cite[3.8.B]{Gro} for discrete groups.
Before stating them, we will make the following \textbf{standing assumption} on a family $\mathcal X$ of  metric spaces:

\noindent{\bf(*)} They must consist of uniformly convex spaces and be stable under the following natural operations:
 ultraproducts; rescaling of the metric;  taking closed geodesically convex subspaces (resp.\ closed vector subspaces when working with Banach spaces).

The main examples one can have in mind are the following classes:
\begin{enumerate}
    \item the class of Hilbert spaces;
    \item the class of closed subspaces of $L^p$-spaces, for a fixed $1<p<\infty$;
     \item more generally: any class of uniformly convex Banach spaces which is stable under ultraproducts;
    \item the class of CAT(0) spaces;
    \item the class of $\mathbb R$-trees.
\end{enumerate}
\begin{thmIntro}[See Theorem \ref{thmInSection:positive-displacement}]\label{thm-positive-displacement}
Let $\mathcal X$ be a family of uniformly convex metric spaces, closed under taking rescaled ultraproducts and closed under taking closed convex subspaces.
Let $G$ be a compactly generated locally compact group such that every continuous isometric action on an element of $\mathcal X$ as almost fixed points. Then $G$ has property F$\mathcal{X}$.
\end{thmIntro}
This theorem was proved by Shalom \cite[Theorem 6.1]{Shalom} for isometric actions of locally compact groups on Hilbert spaces (Shalom states it in terms of reduced cohomology).
\begin{thmIntro}[See Theorem \ref{thmInSection:FCAT(0)compactlypresented}]\label{thm:FCAT(0)compactlypresented}
Let $\mathcal X$ be a family of uniformly convex metric spaces, closed under taking rescaled ultraproducts and closed under taking closed convex subspaces.
Let $G$ be a compactly generated locally compact group with property F$\mathcal{X}$. Then there exists a compactly presented locally compact group $G'$ with property F$\mathcal{X}$, and a continuous surjective homomorphism $G'\to G$ with discrete kernel.
\end{thmIntro}

Two remarks are in order:
\begin{itemize}
    \item Theorems \ref{thm-positive-displacement} and \ref{thm:FCAT(0)compactlypresented} have natural versions for bornological groups: see Corollary \ref{cor-ultraproduct-groups-without-FH} and Theorem \ref{thm-boundedly-presented}.
    \item Theorem \ref{thm:FCAT(0)compactlypresented} is a special case of a more general result saying that property F$\mathcal{X}$ is open in the relevant \emph{marked group topology} for locally compact compactly generated groups (see Theorem \ref{thm:Topen}). 
\end{itemize}
As a special case, we recover the following known result. 
\begin{corollary}\label{cor:Tcompactlypresented}
Let $G$ be a locally compact group with property (T). Then there exists a compactly presented locally compact group $G'$ with property (T), and a continuous surjective homomorphism $G'\to G$ with discrete kernel.
\end{corollary}

Corollary \ref{cor:Tcompactlypresented} has been proved  by Shalom \cite[Theorem 6.7]{Shalom} for discrete groups. It has later been proved by Fisher and Margulis for locally compact groups \cite[Theorem 2.4.]{FM}.  Recently, de la Salle proved a version of Theorem \ref{thm:FCAT(0)compactlypresented} for isometric actions of locally compact groups on certain families of Banach spaces \cite[Corollary 5.13.]{dlS}. Like ours, the previous proofs are based on an ultralimit argument. A difficulty  comes from the fact that when defined naïvely, an ultralimit of isometric actions may not be continuous. However, it remains easy to ensure that orbits of compact sets are bounded, so in our case, we simply can apply Theorem \ref{thm:topo/coarseFX} to replace it by a continuous one. The main interest of our proof is that it applies to isometric actions on non-linear geodesic metric spaces.

\subsection{Bornological groups with bounded geometry}\label{secIntro:BoundedGeom}
We now consider a property of bornological groups which can be seen a coarse analogue of being locally compact for topological groups. 
\begin{definition}
A bornological group $G$ has \emph{bounded geometry} if there is a bounded set $C$ such that all bounded sets $A$ can be covered by finitely many left translates of $C$.
Such a $C$ is called \emph{gordo} if it is also symmetric.
\end{definition}

\begin{remark}
If $C$ is gordo, then all bounded sets $A$ can also be covered by finitely many right translates of $C$, for if $A^{-1}\subseteq g_1C\cup\cdots\cup g_kC$ then $A\subseteq C^{-1}g_1^{-1}\cup\cdots\cup C^{-1}g_k^{-1}$.
\end{remark}

 \begin{remark}
If a bornological group $G$ is such that its bornology is generated by a bounded set $S$, then its bounded structure coincides with $\mathcal B_\sigma$ where $\sigma$ is the action of $G$ on its Cayley graph Cay$(G,S)$. In other words, its coarse structure is the one associated with the left-invariant word metric associated to $S$. Moreover $G$ has bounded geometry if and only if $S^n$ is contained in a union of finitely many left-translates of $S$ for all $n\geq 1$. 
 \end{remark}
 
  Equipped with the bornology $\mathcal K$, a locally compact group has bounded geometry: a gordo set is given by any compact neighbourhood of the identity. More generally, any subgroup $H$ of a locally compact group $G$, equipped with the induced bornology from $G$, has bounded geometry. 
 
A more interesting source of examples are Polish groups with bounded geometry in the sense of Rosendal \cite[Chapter 5]{Ros21}. Lots of explicit examples can be obtained from \cite[Theorem 3]{MR}, and we will see one of them shortly. 

Recall that by a well-known theorem of Serre,
if we have short exact sequence of locally compact groups
\[1\to C\to G\to Q\to 1,\]
such that $Q$ has (T), and $C$ is central and contained in $\overline{[G,G]}$, then $G$ has (T) \cite[Theorem 1.7.11]{BHV07}. A more general version of this result holds for property F$\mathcal E$ where $\mathcal E$ is any class of uniformly convex Banach space which is stable under ultraproduct \cite[Theorem 6.3.]{DT}.  
We do not know whether \cite[Theorem 6.3.]{DT} holds for all bornological groups. However, we are able to prove a version of Serre's theorem for bornological groups with bounded geometry. 
For $A,B$ subsets of a coarse group $G$, we say that $A$ is \emph{coarsely contained} in $B$ if there is a bounded $C\subseteq G$ with $A\subseteq C\cdot B$.

\begin{thmIntro}[See Theorem \ref{thmInSection:Serre}]\label{thm:Serre}
Let $(G,\mathcal B)$ be a bornological group with bounded geometry, whose bornology is generated by a single set.
Let $p\colon G\to Q$ be a surjective morphism and equip $Q$ with the induced bornology $p(\mathcal B)$.
Assume that the kernel $A=\ker(p)$ is central in $G$ and coarsely contained in the commutator $[G,G]$.
Then $G$ has coarse property (T) if and only if $Q$ has coarse property (T).
\end{thmIntro}
Note that even for locally compact groups equipped with $\mathcal K$, the condition on $A$ is weaker than the one in Serre's theorem, so in particular we recover it.

To showcase the power of Theorem \ref{thm:Serre}, we will use it to study certain Polish groups considered in \cite{Ros21}. First, let $\Homeo^+_\mathbb Z(\mathbb R)$ be the group of all lifts of orientation-preserving homeomorphisms $h$ of the circle $S^1$, equipped with the compact open topology.
It is one explicit example of a Polish group with bounded geometry that can be obtained from \cite[Theorem 3]{MR}.
Equipped with $\OB$, this Polish group is coarsely equivalent to $\Z$, and hence has bounded geometry (see \cite[\S3.4.]{Ros21}). 
 \begin{thmIntro}[See Corollary \ref{corInSection:HomeoZR}]\label{thm:HomeoZR}
 The group $\Homeo^+_{\mathbb Z}(\mathbb R)$ equipped with $\OB$ has coarse property (T). Hence, it also has coarse property FH (and in particular topological property FH).
\end{thmIntro}
We get similar results for groups such as the group $\Aut_\mathbb Z(\mathbb Q)$ of lifts of order-preserving permutations of $\mathbb Q/\Z$ and the group $AC^*_\mathbb Z(\mathbb R)$ of lifts of absolutely continuous orientation-preserving homeomorphisms of the circle whose inverses are absolutely continuous.

We turn to an important step in the proof of Theorem \ref{thm:Serre}.
Let us say that an isometric action of a bornological group is \emph{well controlled} if the orbits of bounded sets are relatively compact. Note that this notion is more demanding than being controlled. In particular it is non-empty for linear actions: if $G$ is a locally compact non-discrete group, then the regular representation of $G$ as a discrete group (on $\ell^2(G)$) is never well controlled. 
\begin{definition}
A bornological group $G$ has \emph{coarse property (T-)} if every well-controlled representation with almost invariant vectors has an invariant vector.
It has \emph{coarse property FH-} if every well-controlled isometric action on a Hilbert space has a fixed point.
\end{definition}
As for coarse property (T), we prove the relevant version of  Delorme-Guichardet's theorem (See Theorem \ref{prop-Delorme-Guichardet}).

Clearly coarse property (T) implies coarse property (T-). Since both properties are stronger than usual property (T) for locally compact groups, we deduce from Theorem  \ref{thm:topo/coarseFX} that they are equivalent for locally compact groups. 
It turns out that this equivalence can be extended to bornological groups with  bounded geometry.
\begin{thmIntro}[See Theorem \ref{thmInSection-(T)-and-(T-)-equivalent}]\label{thm-(T)-and-(T-)-equivalent}
For bornological groups with bounded geometry, (T) and (T-) are equivalent.
Similarly, FH and FH- are equivalent for these groups.
\end{thmIntro}
\begin{proposition}[See Proposition \ref{propInSection:wellcontrolledSubaction}]
Let $G$ be a bornological group with bounded geometry and let $\alpha \colon G\curvearrowright \H$ be a controlled affine action on a Hilbert space.
There is a non-empty invariant closed subspace $\H'$ such that the restriction of $\alpha$ to $\H'$ is well controlled.
\end{proposition}\label{prop:wellcontrolledSubaction}

They are also equivalent for the bounded products appearing in Theorem \ref{thm-bounded-products-coarse-T}.
In fact, we have no example where they are not equivalent, although we suspect they exist.

Theorem \ref{thm-(T)-and-(T-)-equivalent} is a consequence of the following result of independent interest, which is in the same spirit as Proposition \ref{prop:continuousSubaction}.

Contrary to Proposition \ref{prop:continuousSubaction}, which has a version for actions on very general geodesic metric spaces, we do not know whether this one holds beyond actions on Hilbert spaces. This comes from the method of proof: Proposition \ref{prop:continuousSubaction} exploits a notion of ``centre'' of bounded subsets which makes sense in a wide class of metric spaces. By contrast, the proof of  Proposition \ref{prop:wellcontrolledSubaction} is based on a more subtle notion of centre: shopping centres, which we could only define in Hilbert spaces.

\subsection{Bounded products}\label{secIntro:BoundedProducts}

We now turn to a class of topological groups for which topological (T) and coarse (T) are distinct (and similarly for property FH). 
Let $\Gamma_n$ be a sequence of (discrete) groups, and let $S_n$ be a symmetric set of generators of $\Gamma_n$.
For each $x\in \Gamma_n$ we define the length $l(x)$ to be the least natural number such that $x\in S_n^{l(x)}$.

\begin{definition}
The \emph{bounded product} of a sequence of groups $\Gamma_n$ with symmetric generating sets $S_n$ is the group
\[\prod_n\left\{\Gamma_n,S_n\right\}  = \{(g_n) \in \prod_n\Gamma_n\mid \sup_nl(g_n) < \infty\}.\]
\end{definition}

Let $G = \prod_n\left\{\Gamma_n,S_n\right\}$.
This has a canonical bornology $\mathcal B$, where $A\subseteq G$ is bounded if and only if $\sup_{(g_n)\in A}\sup_nl(g_n) < \infty$. 

There is also a natural topology on $G$: for any proper function $\rho\colon\mathbb N\to \mathbb Z_{\geq 0}$, define
\[U_\rho = \{(g_n)\in G \mid l(g_n) \leq \rho(n)\}.\]
The sets $U_\rho$, together with their translates, form a basis for the topology on $G$.
\begin{proposition}[See Proposition \ref{propInSection-bornologies-bounded-product}]
\label{prop-bornologies-bounded-product}
With respect to this topology, we have $\mathcal B = \mathcal K = \OB$, provided that all $S_n$ are finite.
\end{proposition}

We shall say that a topological group $G$ has {\it topological} property (T) if all continuous unitary representations with almost invariant vectors (relative to compact subsets) have invariant vectors (see Definition \ref{def-(T)-for-topological-bornological-group}).

If the $\Gamma_n$ are all finite groups and $|S_n|$ is bounded, we can characterise exactly when $G$  has coarse or topological property (T).
\begin{thmIntro}[See Theorem \ref{thmInSection-bounded-products-topological-T}]\label{thm-bounded-products-topological-T}
Let $\Gamma_n$ be a sequence of finite groups with symmetric generator sets $S_n$.
Suppose that $\sup_n|S_n|<\infty$.
Let $X$ be the coarse disjoint union of the Cayley graphs $\Cay(\Gamma_n,S_n)$.
Then the bounded product $G$ has topological property (T) if and only if $X$ is an expander sequence.
\end{thmIntro}

Given a non-principal ultrafilter $\mathcal U$ on $\mathbb N$ there is an equivalence relation $\sim$ on $G$, with $(g_n)\sim(h_n)$ if $\{n\mid g_n=h_n\}\in\mathcal U$.
Taking the quotient by this equivalence relation gives a group $G_\mathcal U = \prod_n\{\Gamma_n,S_n\}/\mathcal U$, generated by the image of $\prod_nS_n$.
We will view it as a discrete group. If $|S_n|=k$ is a fixed integer, then the groups obtained as $G_\mathcal U$ have a natural interpretation in terms of the space of marked groups, i.e.\ quotients of the free groups on $k$ generators: indeed, it consists of the closure of $\{(\Gamma,S_n)\mid n\in\mathbb N\}.$

\begin{thmIntro}[See Theorem \ref{thmInSection-bounded-products-coarse-T}]
\label{thm-bounded-products-coarse-T}
Let $\Gamma_n$ be a sequence of finite groups with symmetric generator sets $S_n$.
Suppose that $\sup_n|S_n|<\infty$.
Let $X$ be the coarse disjoint union of the Cayley graphs $\Cay(\Gamma_n,S_n)$.
The following are equivalent:
\begin{enumerate}
    \item The coarse space $X$ has geometric property (T).
    \item All ultraproducts $G_\mathcal U = \prod_n\Gamma_n/\mathcal U$ have property (T) with uniform spectral gap.
    \item All ultraproducts $G_\mathcal U$ have property (T).
    \item There is a finitely generated group $\Gamma$ with property (T) such that all $\Gamma_n$ are a quotient of $\Gamma$.
    \item The bounded product $G$ has coarse property (T).
    \item The bounded product $G$ has coarse property (T-).
\end{enumerate}
\end{thmIntro}

The equivalence between (i) and (ii) are proved in \cite[Theorem 1]{MOSS}.
To the best of our knowledge the equivalence between (i) and (iv) is new, but is reminiscent of a similar statement proved in \cite{WY}: if $\Gamma_n=\Gamma/\Lambda_n$, where $\Lambda_n$ is a decreasing sequence of finite index normal subgroups of $\Gamma$ such that $\cap_n\Lambda_n=\{1\}$, then the coarse space $X$ has geometric property (T) if and only if $\Gamma$ has (T).

\subsection*{Acknowledgement}
We would like to thank Tim de Laat for his valuable remarks on a preliminary version of the paper. 
The second author is supported by the Deutsche Forschungsgemeinschaft under Germany's Excellence Strategy - EXC 2044 - 390685587, Mathematics M\"unster: Dynamics - Geometry - Structure.

\section{Centres in uniformly convex metric spaces}\label{sec:centres}
Let $E$ be a Banach space.
It is \emph{reflexive} if the natural map to the double dual $E\to E^{**}$ is an isomorphism.
It is called superreflexive if all its ultrapowers are reflexive.
A Banach space is superreflexive if and only if it has an equivalent norm that is uniformly convex, i.e. for every $\epsilon>0$ there is $\delta > 0$ such that if $x,y$ are unit vectors with $\norm{x-y} \geq \epsilon$, then $\norm{\frac{x+y}2} < 1-\delta$. Moreover, for any isometric action $\sigma$ of a group $G$ on a  superreflexive space $X$, there exists an equivalent uniformly convex norm $\|\cdot\|'$ on $X$ such that $\sigma$ is isometric on $(X,\|\cdot\|')$, see \cite{Enflo}.
Examples of superreflexive Banach spaces are Hilbert spaces and all $L^p$ spaces with $1< p < \infty$.

Even more generally, we can consider a uniquely geodesic complete metric space $X$.
For any $x,y\in X$, the geodesic between them is denoted by $[x,y]$ and the midpoint by $m_{xy}$.
The distance between $x$ and $y$ is denoted by $|xy|$.
A subset $A\subseteq X$ is \emph{convex} if $x,y\in A$ implies $[x,y]\subseteq A$.
For any $A\subseteq X$, the \emph{convex hull} $\co(A)$ is the smallest convex set that contains it.
Its closure is denoted by $\cco(A)$.

We take the following notion of uniform convexity from \cite[Definition 7.1]{Nic08}.
\begin{definition}\label{def:uniformly convex}
A uniquely geodesic complete metric space $X$ is \emph{uniformly convex} if for each $\epsilon>0$ there is $\delta >0$ such that for all $x,y,z \in X$ with $|xy| \geq\epsilon\max(|xz|,|yz|)$, we have $|zm_{xy}|\leq (1-\delta)\max(|xz|,|yz|)$.
\end{definition}
We will say ``uniformly convex metric space" as short-hand for ``uniformly convex uniquely geodesic complete metric space".
This generalizes uniform convexity for Banach spaces.
As a matter of notation, we will write $\delta_\epsilon$ for a $\delta$ that works.

\subsection{Centres of bounded sets}
Let $X$ be a uniformly convex metric space.
Let $A\subseteq X$ be a non-empty bounded subset.
For $x\in X$ let $r(x) = \sup_{a\in A}|ax|$.
Let $\rho(A) = \inf_{x\in X}r(x)$.
As shown in \cite[Proposition 7.2]{Nic08}, there is a unique point $Z(A)$, called the centre, that satisfies $r(Z(A)) = \rho(A)$.
There are no additional regularity requirements for $A$.
It is an isometric invariant, meaning that whenever $u\in \Aut(X)$ is an isometry, we have $Z(uA) = uZ(A)$.
In this section we describe this centre and some of its variants, which will be the key ingredient in many of our proofs.

For convenience we provide the construction of the centre.

\begin{lemma}\label{lem-of-centre}
Let $A\subseteq X$ be a non-empty bounded subset of a uniformly convex metric space.
There is a unique point $Z(A)$ satisfying $r(Z(A)) = \rho(A)$.
\end{lemma}
\begin{proof}
Let $x_n\in X$ be a sequence of points satisfying $\lim_{n\to\infty}r(x_n) = \rho(A)$.
We will show that it is a Cauchy sequence.
Let $\epsilon>0$ and let $m,n$ large enough that $r(x_m),r(x_n)\leq (1+\delta_\epsilon)\rho(A) \leq 2\rho(A)$.
Suppose that $|x_mx_n|\geq 2\epsilon\rho(A)$.
Let $y$ be the midpoint between $x_m$ and $x_n$, then we get by uniform convexity $|ya| \leq (1-\delta_\epsilon)\max(|x_ma|,|x_na|) \leq (1-\delta_\epsilon^2)\rho(A)$ for all $a\in A$.
Thus $r(y) < \rho(A)$, a contradiction.
So we have $|x_mx_n| < 2 \epsilon \rho(A)$, showing that $x_n$ is a Cauchy sequence.
Then the limit $x$ satisfies $r(x) = \rho(A)$, and the same inequality shows that $x$ is unique with this property.
\end{proof}

\begin{lemma}\label{lem-additional-information-centre}
Let $A\subseteq X$ be a bounded subset of a uniformly convex metric space.
\begin{enumerate}[(i)]
    \item Let $\epsilon>0$ and consider the annulus $D = \{a\in A\mid |a,Z(A)| \geq \rho(A)-\epsilon\}$.
    We have $Z(D) = Z(A)$.
    \item Let $B\subseteq A$ non-empty and $\epsilon > 0 $.
    Then
    \[ |Z(A)Z(B)| \leq \epsilon \rho(A) + \kappa_\epsilon(\rho(A)-\rho(B))\]
    where $\kappa_\epsilon = \frac{2-\epsilon}{1-\exp\left(\frac{-\log(\epsilon/2)\log(1-\delta_\epsilon)}{\log2}\right)}$ depends only on $\epsilon$ and $\delta_\epsilon$.
    \item If $X$ is a Hilbert space then $Z(A) \in \cco(A)$.
\end{enumerate}
\end{lemma}
\begin{proof}
\begin{enumerate}
    \item Since $D\subseteq A$ we have $\rho(D)\leq \rho(A)$, and we also have $\sup_{d\in D}|dZ(A)|= \rho(A)$.
    Assume for a contradiction that $Z(D) \neq Z(A)$.
    Let $x$ be a point on the geodesic $[Z(D),Z(A)]$ with $0< |xZ(A)| < \epsilon$.
    Let $m$ be the midpoint between $x$ and $Z(A)$.
    Let $\epsilon' = \frac{|xZ(A)|}{\rho(A)}$.
    For $a\in D$ we have $|aZ(D)| \leq \rho(D)\leq \rho(A)$ and $|aZ(A)| \leq \rho(A)$.
    By convexity, $|ax| \leq \rho(A)$.
    By uniform convexity, $|am| \leq (1-\delta_{\epsilon'})\rho(A)$.
    For $a\in A\setminus D$ we have $|aZ(A)| \leq \rho(A) - \epsilon$ so $|am| \leq \rho(A) - \frac\epsilon2$.
    This shows that $m$ is a better centre for $A$ then $Z(A)$, giving a contradiction.
\item Let $\epsilon > 0$ and suppose $|Z(A)Z(B)| \geq \epsilon\rho(A) \left(\frac{\rho(A)}{\rho(B)}\right)^\frac{\log(2)}{-\log(1-\delta_\epsilon)}$.
Let $n = \lfloor \frac{\log(\rho(A)) - \log(\rho(B))}{-\log(1-\delta_\epsilon)} + 1 \rfloor$.
Then $(1-\delta_\epsilon)^n < \frac{\rho(A)}{\rho(B)}$, and $|Z(A)Z(B)| \geq 2^{n-1}\epsilon\rho(A)$.
Let $m_0 = Z(A)$ and for $1\leq k\leq n$, let $m_k$ be the midpoint between $m_{k-1}$ and $Z(B)$.
For $b\in B$ we have $|Z(B)b| \leq \rho(B)$ and $|Z(A)b|\leq \rho(A)$.
By uniform convexity and induction we get $|m_kb| \leq \max((1-\delta_\epsilon)^k\rho(A), (1-\delta_\epsilon)\rho(B))$.
We get $\sup_{b\in B}|m_nb| \leq \max((1-\delta_\epsilon)^n\rho(A) , (1-\delta_\epsilon)\rho(B)) < \rho(B)$.
Thus $m_n$ is a better centre for $B$ then $Z(B)$, a contradiction.

So $|Z(A)Z(B)| < \epsilon\rho(A)\left(\frac{\rho(A)}{\rho(B)}\right)^\frac{\log(2)}{-\log(1-\delta_\epsilon)}$.
Since any point in $B$ is at distance at most $\rho(A)$ from both $Z(A)$ and $Z(B)$ we also have $|Z(A)Z(B)|\leq 2\rho(A)$.
Using the inequality $\max\left(\left(\frac1{1-x}\right)^\frac{\log(2)}{-\log(1-\delta_\epsilon)},2\epsilon^{-1}\right)\leq 1+\epsilon^{-1}\kappa_\epsilon x$ for $x\in[0,1)$ we get $|Z(A)Z(B)|\leq \epsilon\rho(A)+\kappa_\epsilon(\rho(A)-\rho(B))$.
\item Suppose $Z(A) \not\in\cco(A)$.
Then there is a unit vector $w$ with $\langle Z(A), w \rangle < \inf_{a\in A}\langle a,w\rangle$.
But then $Z(A) + \epsilon w$ is a better centre for $A$ for small $\epsilon$, contradiction.
\end{enumerate}
\end{proof}
\begin{remark}
If $X$ is a Hilbert space, we can find the simpler estimate $|Z(A)Z(B)| \leq \sqrt{\rho(A)^2-\rho(B)^2}$ if $B\subseteq A$.
\end{remark}

 We record the following easy fact.
\begin{proposition}\label{prop-dense-subgroup}
Let $G$ be a bornological group and $H\subseteq G$ a coarsely dense subgroup. Let $G$ act isometrically and controlled on a uniformly convex metric space $X$ and suppose there is a point in $X$ fixed by $H$.
Then there is also a point fixed by $G$.
\end{proposition}
\begin{proof}
There is a bounded $C\subseteq G$ with $CH = G$.
Let $\alpha\colon G\curvearrowright X$ be the isometric action and $x\in X$ the point which is fixed by $H$.
Now $\alpha(G)x = \alpha(C)x$ which is bounded, so $\alpha$ has bounded orbits. 
By Lemma \ref{lem-of-centre}, there is a unique centre $Z(\alpha(G)x)$.
This is a fixed point for $G$: of any $g\in G$ we have $gZ(\alpha(G)x) = Z(g\alpha(G)x) = Z(\alpha(G)x)$.
\end{proof}

\subsection{Mean centre}
The next variant of the centre allows us to talk about some kind of average of a mean on a uniformly convex Banach space, without additional regularity properties.
It generalises both the centre of a bounded set, and the barycentre of a probability measure. This notion will be used to provide an alternative proof of Proposition \ref{prop:continuousSubaction} in the special case of an isometric action of a locally compact group on a uniformly convex Banach space. However we believe it could be useful for future applications. 

\begin{construction}\label{constr:meancentre}
Let $X$ be a closed convex subset of a uniformly convex Banach space.
Let $C\subseteq X$ be a bounded subset.
Let $\Sigma$ be an algebra on $C$ and let $\mu\colon\Sigma\to [0,\infty)$ be a finitely additive measure with $\mu(C) > 0$.
We construct the \emph{mean centre} $\tilde Z(\mu) \in X$ with the following properties:
\begin{enumerate}[(i)]
    \item for $A\in \Sigma$ we have
\[ \tilde Z(\mu) = \frac{\mu(A)}{\mu(C)}\tilde Z(\mu_{|A}) + \frac{\mu(C\setminus A)}{\mu(C)}\tilde Z(\mu_{|C\setminus A});\]
\item For any isometry $g\in \Isom(X)$ we have $\tilde Z(g_\ast\mu)=g\cdot \tilde Z(\mu)$.
\end{enumerate}
\end{construction}
\begin{proof}
For a bounded subset $B\subseteq X$, recall the definitions of $Z(B)$ and $\rho(B)$ from Lemma \ref{lem-of-centre}.
For any finite partition $P = (A_i)$ of $C$ into subsets in $\Sigma$, we define
\[\tilde Z_P(\mu) = \frac1{\mu(C)}\sum_i\mu(A_i)Z(A_i)\]
and
\[\overline\rho_P(\mu) = \frac1{\mu(C)}\sum_i\mu(A_i)\rho(A_i).\]
We will define $\tilde Z(\mu)$ as a limit of $\tilde Z_P(\mu)$, so we have to show convergence as $P$ gets more refined.

Let $\epsilon > 0$ and $\kappa_\epsilon$ as in Lemma \ref{lem-additional-information-centre}(ii).
Then $\norm{Z(A)-Z(B)} \leq \epsilon\rho(A) + \kappa_\epsilon(\rho(A)-\rho(B))$ for all $\emptyset\neq B\subseteq A\subseteq C$.

Let $P = (A_i)$ be a finite partition of $C$ and let $Q = (B_{ij})$ be a refinement of $P$.
Then we have
\begin{align*}
    \norm{\tilde Z_P(\mu)-\tilde Z_Q(\mu)}&= \frac1{\mu(C)}\norm{\sum_i\left(\sum_j\mu(B_{ij})Z(B_{ij})\right)-\mu(A_i)Z(A_i)}\\
    &\leq \frac1{\mu(C)}\sum_{i,j}\mu(B_{ij})\norm{Z(B_{ij}) - Z(A_i)}\\
    &\leq \frac{\epsilon}{\mu(C)}\sum_{i,j}\mu(B_{ij})(\rho(A_i)+\epsilon^{-1}\kappa_\epsilon(\rho(A_i)-\rho(B_{ij})))\\
    & = \epsilon\overline\rho_P(\mu)+\kappa_\epsilon(\overline \rho_P(\mu) - \overline \rho_Q(\mu)).
\end{align*}
Besides, we observe that $\overline\rho_Q(\mu) \leq\overline\rho_P(\mu)$.
Now let $P_1,P_2,\ldots$ be a sequence of partitions, such that $P_{n+1}$ is a refinement of $P_n$, and $\overline\rho_{P_n}(\mu)$ converges to the infimum of $\overline\rho_P(\mu)$ taken over all partitions $P$.
Then the computation above shows that $\tilde Z_{P_n}(\mu)$ is a Cauchy sequence: for large $m,n$, we will have $\overline\rho_{P_n}(\mu) - \overline\rho_{Q_n}(\mu) < \frac\epsilon{\kappa_\epsilon}$, and then $\norm{\tilde Z_{P_n}(\mu) - \tilde Z_{P_m}(\mu)} < \epsilon\rho(C) + \epsilon$, and this works for any $\epsilon > 0$.
So the sequence $\tilde \rho_{P_n}(\mu)$ converges to a point $\tilde Z(\mu) \in X$ which is independent of the chosen sequence $P_n$.

Now let $A\in \Sigma$.
We can choose the partitions $P_n$ to be refinements of the partition $\{A,U\setminus A\}$.
Then they split into a partition of $A$ and one of $U\setminus A$, and we see that
\[\tilde Z(\mu) = \frac{\mu(A)}{\mu(U)}\tilde Z(\mu_{|A}) + \frac{\mu(U\setminus A)}{\mu(U)}\tilde Z(\mu_{|U\setminus A}).\]

The equivariance of $\tilde Z$ by isometries follows directly from the definition.
\end{proof}

\begin{remark}
If $\Sigma=\{\emptyset, C\}$, then $\tilde Z(\mu)=Z(C)$. Moreover, one can check that if $C$ is weakly compact, and $\mu$ is a Radon probability measure on $C$ (equipped with the weak topology), then  $\tilde Z(\mu)$ coincides with the barycentre of $\mu$. 
\end{remark}

\begin{remark}
The above construction does not quite work for any uniformly convex metric space $X$ because there is no linear structure on $X$. However, instead of constructing $\tilde Z(\mu)$ as a point in $X$, it is possible by a similar construction to construct a Radon probability measure on $X$.
\end{remark}

\subsection{Shopping centres}\label{sec:shoppingcentres}

In this section, we prove a technical result which will be central to the proof of Proposition \ref{prop:wellcontrolledSubaction}, hence of Theorem \ref{thm-(T)-and-(T-)-equivalent}. Recall that given a controlled action $\sigma$ of a bornological group $G$ with bounded geometry on a Hilbert space $\mathcal H$, we need to find $v\in \mathcal H$ such that $\sigma(A)\cdot v$ is relatively compact for all bounded $A\subset G$. Example \ref{ex:Znot enough} below gives a hint that simply taking the centre $Z(\sigma(A))$ for some gordo subset $A$ does not seem to be a good approach. Instead, we will need another concept of centre which tries to ignore finite-dimensional subspaces of $\H$. 

For a finite-dimensional subspace $V\subseteq \H$, let $p_{V^\perp}$ be the projection onto $V^\perp$.
For a bounded set $\Lambda \subseteq \H$, let $s_V(\Lambda) = \rho(p_{V^\perp}(\Lambda))$.
Let $s(\Lambda) = \inf_Vs_V(\Lambda)$, taking the infimum over all finite-dimensional subspaces. Note that $s(\Lambda)=0$ if and only if $\Lambda$ is relatively compact.

\begin{definition}
Let $\Lambda\subseteq \H$ be a bounded subset.
A point $z\in \H$ is called a \emph{shopping centre} for $\Lambda$ if for all finite-dimensional $V\subseteq \H$ and all $\epsilon>0$ we have $z\in\cco(\{v\in\Lambda\mid\norm{p_{V^\perp}(v-z)}\geq s(\Lambda)-\epsilon\})$.
\end{definition}

We first prove the existence of shopping centres.

\begin{lemma}\label{lem-existence-shopping-centres}
Every non-empty bounded subset of a Hilbert space has a shopping centre.
\end{lemma}
\begin{proof}
Let $\Lambda \subseteq \H$ be bounded.
We write $s_V = s_V(\Lambda)$ and $s = s(\Lambda)$.
We can assume that $s > 0$, because the statement is trivial otherwise.
For any finite-dimensional $V\subseteq \H$ and $\epsilon>0$ define 
\[ A_{V,\epsilon} = \{a\in\H \mid a\in \cco(\{v\in\Lambda\mid \norm{p_{V^\perp}(a-v)} \geq s-\epsilon\})\}.\]
First we show that $A_{V,\epsilon}$ is non-empty.

Let $x = Z(p_{V^\perp}(\Lambda)) \in V^\perp$.
Since $\rho(p_{V^\perp}(\Lambda)) = s_V \geq s$, we have
\[x \in \cco(\{p_{V^\perp}v\mid v\in\Lambda,\ \norm{p_{V^\perp}v-x}\geq s-\tfrac\epsilon2\})\]
by Lemma \ref{lem-additional-information-centre}(i) and (iii).
So we can find $y \in \co(\{p_{V^\perp}v\mid v\in\Lambda,\ \norm{p_{V^\perp}v-x} \geq s-\frac\epsilon2\})$ with $\norm{y-x}<\frac\epsilon2$.
Write $y = \sum_{i=1}^kc_ip_{V^\perp}v_i$ with $c_i\geq 0$ and $\sum_{i=1}^kc_i = 1$ and $v_i\in\Lambda$ and $\norm{p_{V^\perp}v_i-x}\geq s-\frac\epsilon2$ for all $i$.
Let $a = \sum_{i=1}^kc_iv_i$.
Then
\[a \in \co(\{v\in\Lambda \mid \norm{p_{V^\perp}v-x} \geq s-\tfrac\epsilon2\}).\]
We have $\norm{p_{V^\perp}a - x} = \norm{y -x} <\frac\epsilon2$, so 
\[\{v\in\Lambda\mid\norm{p_{V^\perp}v-x}\geq s-\tfrac\epsilon2\}\subseteq \{v\in\Lambda\mid\norm{p_{V^\perp}(v-a)}\geq s-\epsilon\},\]
so $a\in A_{V,\epsilon}$.

Note that, if $\epsilon'\leq\epsilon$ and $V'\supseteq V$, we have $A_{V',\epsilon'}\subseteq A_{V,\epsilon}$.
So now we already know that any finite intersection of the $A_{V,\epsilon}$ is non-empty.

We can bound the diameter of $p_{V^\perp}A_{V,\epsilon}$.
Assume $\epsilon < s$.
Let $x = Z(p_{V^\perp}\Lambda) \in V^\perp$ and let $a\in A_{V,\epsilon}$.
Since $a\in \cco(\{v\in\Lambda\mid \norm{p_{V^\perp}(a-v)} \geq s-\epsilon\})$, we can find, for every $\eta>0$, a vector $v\in\Lambda$ with $\norm{p_{V^\perp}(a-v)} \geq s-\epsilon$ and $\langle p_{V^\perp}(a-v),p_{V^\perp}(a)-x\rangle \leq \eta$.
Then $\norm{x-p_{V^\perp}(v)}^2 = \norm{p_{V^\perp}(a-v)}^2 + \norm{x-p_{V^\perp}(a)}^2-2\langle p_{V^\perp}(a-v),p_{V^\perp}(a)-x\rangle\geq (s-\epsilon)^2 + \norm{x-p_{V^\perp}(a)}^2 - 2\eta$.
On the other hand we have $\norm{x-p_{V^\perp}(v)}^2\leq s_V^2$, so $\norm{p_{V^\perp}(a)-x}^2 \leq s_V^2 - (s-\epsilon)^2 + 2\eta$.
This holds for all $\eta > 0$ so $\norm{p_{V^\perp}(a)-x} \leq \sqrt{s_V^2 - (s-\epsilon)^2}$.
So the diameter of $p_{V^\perp}A_{V,\epsilon}$ is at most $2\sqrt{s_V^2 - (s-\epsilon)^2}$.
This then also holds for the diameter of the weak closure of $p_{V^\perp}A_{V,\epsilon}$.

Let $b$ be in the weak closure of $A_{V,\epsilon}$.
Then for any $a\in A_{V,\epsilon}$ we have $\norm{p_{V^\perp}(a-b)} \leq 2\sqrt{s_V^2 - (s-\epsilon)^2}$, so
\[ \{ v\in \Lambda \mid \norm{p_{V^\perp}(a-v)} \geq s-\epsilon\} \subseteq \{v\in\Lambda\mid \norm{p_{V^\perp}(b-v)}\geq s-\epsilon - 2\sqrt{s_V^2-(s-\epsilon)^2}\}.\]
This shows that
\[A_{V,\epsilon} \subseteq \cco(\{v\in\Lambda \mid\norm{p_{V^\perp}(b-v)}\geq s-\epsilon - 2\sqrt{s_V^2-(s-\epsilon)^2}\}).\]
Since convex closed sets are also weakly closed, we find that $b\in A_{V,\epsilon + 2\sqrt{s_V^2-(s-\epsilon)^2}}$.
So the weak closure of $A_{V,\epsilon}$ is contained in $A_{V,\epsilon + 2\sqrt{s_V^2-(s-\epsilon)^2}}$.

This shows that any $A_{V,\epsilon}$ also contains the weak closure of another $A_{V',\epsilon'}$: just take $\epsilon'$ small enough that $\epsilon' + 2\sqrt{(s+\epsilon')^2 - (s-\epsilon')^2} \leq \epsilon$ and $V'$ large enough that $s_{V'}\leq s + \epsilon'$.
Now the weak closure of $A_{V,\epsilon}$ is weakly compact, and any finite intersection of them is non-empty, so there is $z$ that is contained in the weak closure of every $A_{V,\epsilon}$.
By the above, $z\in A_{V,\epsilon}$ for all finite-dimensional $V\subseteq \H$ and $\epsilon > 0$, so $z$ is a shopping centre for $\Lambda$.
\end{proof}

Shopping centres are not unique, but the set of shopping centres is relatively compact.
In fact we can prove a stronger result, which we will use later to prove Proposition \ref{prop:wellcontrolledSubaction}. 

\begin{lemma}\label{lem-mall-relatively-compact}
Let $\Lambda \subseteq \H$ be a non-empty bounded subset.
Let $f_1,\ldots,f_k \colon \H\to \H$ be isometries, and let $\Delta = \bigcup_{i=1}^kf_i(\Lambda)$.
We say $z\in \H$ is a \emph{possible shopping centre} if there exists an isometry $f\colon \H\to \H$ such that $f(\Lambda)\subseteq \Delta$ and $z$ is a shopping centre of $f(\Lambda)$.
Let $M$ be the set of possible shopping centres.
Then $M$ is relatively compact.
\end{lemma}
\begin{proof}
We write $s_V = s_V(\Lambda)$ and $s = s(\Lambda)$.
Let $\epsilon>0$.
Choose a finite-dimensional affine subspace $V\subseteq \H$ with $s_V \leq s+\epsilon$ and $Z(p_{V^\perp}(\Lambda)) \in V$.
Then for all $v\in \Lambda$ we have $d(v,V) \leq s+\epsilon$.

Let $W$ be the affine hull of the union $\bigcup_{i=1}^kf_i(V)$.
This is again a finite-dimensional affine subspace.
For all $v\in f_i(\Lambda)$, we have $d(v,f_i(V)) \leq s+\epsilon$, so for all $v\in\Delta$ we have $d(v,W)\leq s+\epsilon$.

Now let $z\in M$ be a possible shopping centre.
We will show that $z$ is close to $W$. Let $w$ be the projection of $z$ on $W$.
To simplify notation, let us assume that $z=0$.
There is an isometry $f\colon\H\to\H$ such that $f(\Lambda)\subseteq \Delta$ and $0$ is a shopping centre of $f(\Lambda)$.
Then there is $v\in f(\Lambda)$ such that $\langle v,w\rangle \leq \epsilon$ and $\norm{p_{W^\perp}(v)} \geq s-\epsilon$.
Then
\begin{align*}
d(v,W)^2 &= \norm{p_{W^\perp}(v-w)}^2 \\
&= \norm{p_{W^\perp}(v)}^2 + \norm{p_{W^\perp}(w)}^2 - 2\langle p_{W^\perp}(v),w\rangle\\
&= \norm{p_{W^\perp}(v)}^2 + \norm{w}^2 - 2\langle v,w\rangle \\
&\geq (s-\epsilon)\max(0,s-\epsilon) + \norm{w}^2 -2\epsilon.
\end{align*}
On the other hand, we have $v\in\Delta$, so $d(v,W) \leq s+\epsilon$.
So we get $d(0,W) = \norm{w} \leq \eta = \sqrt{4\epsilon\max(s,\epsilon) + 2\epsilon}$.

So $M$ is contained in an $\eta$-neighbourhood of $W$.
Since $M$ is also bounded, we can cover $M$ with finitely many balls of radius $2\eta$.
As $\epsilon$ gets smaller, $\eta$ gets arbitrarily small, so $M$ is totally bounded, which is the same as relatively compact.
\end{proof}

To justify our relatively complicated notion of shopping centre, we give an example showing that Lemma \ref{lem-mall-relatively-compact} is not true for the usual centre.
\begin{example}\label{ex:Znot enough}
Let $\H$ be a separable Hilbert space with basis $v_1,v_2,v_3,\ldots$.
Let $B(0,1)$ denote the closed ball or radius 1 around the origin and let $\Lambda$ be the bounded set $B(0,1)\cup\{10v_1+v_2,-10v_1+v_2\}$.
Let $f_1\colon\H\to \H$ be the identity, let $f_2$ be the isometry given by adding $10v_1$ and let $f_3$ be the isometry given by subtracting $10v_1$.
Then $\Delta = B(0,1)\cup B(10v_1,1)\cup B(-10v_1,1)\cup\{20v_1+v_2,-20v_1+v_2\}$.
Note that $Z(\Lambda) = v_2$.
Now for $n\geq 2$ let $g_n\colon \H\to \H$ be the isometry switching $v_2$ and $v_n$ while fixing all other basis vectors.
Then $g_n(\Lambda) \subseteq \Delta$.
The centre of $g_n(\Lambda)$ is $v_n$.
The set $\{v_2,v_3,...\}$ is not relatively compact, showing that Lemma \ref{lem-mall-relatively-compact} does not hold for the usual centre.
Note that the origin is the unique shopping centre of $\Lambda$.
\end{example}

\section{Topological versus controlled actions}\label{sec-topological-groups}

In this section we will mainly be concerned by topological groups $(G,\tau)$, i.e. groups equipped with a topology making the multiplication and inverse maps continuous.
We will also always assume that the topology is Hausdorff. Each topological group $(G,\tau)$ has two natural bornologies associated to it, $\mathcal K$ and $\OB$, as described in \cite{Ros21}.
\begin{definition}
Let $(G,\tau)$ be a topological group.
The bornology $\mathcal K$ consists of all relatively compact subsets of $G$.
\end{definition}
This is a bornology because the product of compact sets is again compact.

Another bornology comes from considering continuous left-invariant pseudometrics on $G$.
The following proposition is \cite[Proposition 2.7]{Ros21}:
\begin{proposition}\label{prop-coarsely-bounded}
Let $(G,\tau)$ be a topological group and $A\subseteq G$.
The following are equivalent:
\begin{enumerate}
    \item For every continuous left-invariant pseudometric $d$ on $G$ we have $\diam_d(A) < \infty$.
    \item For every continuous isometric action $\alpha\colon G\curvearrowright X$ on a metric space $X$, and every $x\in X$, we have $\diam(\alpha(A)x) < \infty$.
    \item For every sequence of open sets $V_1\subseteq V_2\subseteq V_3\cdots \subseteq G$ of open subsets with $V_n^2\subseteq V_{n+1}$ and $\bigcup_nV_n = G$, we have $A\subseteq V_n$ for some $n$.
\end{enumerate}
\end{proposition}
\qed
\begin{definition}\label{def:OB} The bornology $\OB$ consists of all subsets $A$ of $G$ that satisfy the conditions of Proposition \ref{prop-coarsely-bounded}.
\end{definition}
We have $\mathcal K\subseteq \OB$.
For $\sigma$-compact locally compact groups, we have $\mathcal K = \OB$ (see \cite[Corollary 2.8]{Ros21}).

We want to study the relationship between topological fixed point properties and coarse fixed point properties. We have already defined in the introduction the notions of topological versus coarse fixed point properties. Since the bornologies we shall consider in this paper will always be contained in $\OB$, all continuous isometric actions are controlled, so that topological FX is always stronger than coarse FX (for any metric space $X$).

Now if we turn to property (T), it is less clear how topological (T) should be defined: indeed, even if we focus on continuous representations, there are more than one way of defining almost fixed points. A solution is then to define topological property (T) for a topological group {\it equipped with a bornology} (contained in $\OB$).

\begin{definition}\label{def-(T)-for-topological-bornological-group}
A \emph{topological bornological group} is a triple $(G,\tau,\mathcal B)$ such that $(G,\tau)$ is a topological group and $(G,\mathcal B)$ is a bornological group, and $\mathcal B$ is contained in $\OB(\tau)$.
Let $\mathcal E$ be a family of Banach space. We say $(G,\tau,\mathcal B)$ has \emph{topological property (T$\mathcal E$)} if every continuous (wrt $\tau$) isometric representation on an element of $\mathcal E$ with almost invariant vectors (wrt $\mathcal B$) has a fixed point.
\end{definition}
For any topological group $(G,\tau)$ there are the natural bornologies $\mathcal K$ and $\OB$, and we can consider topological property (T) for the triples $(G,\tau,\mathcal K)$ and $(G,\tau,\OB)$ separately.
For a $\sigma$-compact locally compact group $(G,\tau)$ we have $\mathcal K = \OB$, and the properties above reduce to their usual versions when we take this bornology.

\subsection{Proof of Proposition \ref{prop:continuousSubaction}}\label{sec:prop:continuousSubaction}

To be able to prove anything about topological bornological groups, we need some relation between the topology and the bornology.
It turns out that all we need is the following.
\begin{definition}
A topological bornological group $(G,\tau,\mathcal B)$ is \emph{locally bounded} if there is a non-empty open bounded set.
\end{definition}

For an action of $G$ on a topological space $X$, let $X_{\cont} = \{x\in X\mid g\to gx \text{ is continuous}\}$ be the continuous part of the action.
\begin{lemma}\label{lem:X-cont-subspace}
Let $(G,\tau)$ be a topological group acting isometrically on a uniformly convex metric space $X$.
Then $X_{\cont}$ is an invariant closed convex subspace of $X$.
If $X$ is a Banach space and the action is linear, $X_{\cont}$ is also a linear subspace.
\end{lemma}
\begin{proof}
Since the action is isometric, the continuous part is invariant and $x\in X_{\cont}$ if and only if $g\to gx$ is continuous at the identity.

Let $x\in\overline X_{\cont}$ and $\epsilon>0$.
There is $y\in X_{\cont}$ with $|xy|<\epsilon$ and there is an open neighbourhood of the identity $U\subseteq G$ such that $|gy,y|<\epsilon$ for all $g\in U$.
Then for $g\in U$ we have $|gx,x|\leq |gx,gy|+|gy,y|+|yx|<3\epsilon$, showing that $x\in X_{\cont}$.

To show that $X_{\cont}$ is convex, it is enough to show that if $x,y\in X_{\cont}$ then the midpoint $m_{xy}$ is in $X_{\cont}$.
Let $\epsilon>0$ and let $U$ be a neighbourhood of the identity such that $g\in U$ implies $|gx,x|<\frac12|xy|\delta_\epsilon$ and $|gy,y|<\frac12|xy|\delta_\epsilon$.
Let $g\in U$.
Since the action is isometric we have $m_{gx,gy} = gm_{xy}$.
We have $|x,gm_{xy}|\leq |x,gx|+|gx,gm_{xy}|\leq \frac12|xy|(1+\delta_\epsilon)$.
Let $n$ be the midpoint between $m_{xy}$ and $gm_{xy}$ and suppose that $|m_{xy},gm_{xy}| \geq \frac12\epsilon|xy|(1+\delta_\epsilon)$.
By uniform convexity we have $|xn|\leq (1-\delta_\epsilon)\frac12|xy|(1+\delta_\epsilon)$ and similarly $|yn|\leq (1-\delta_\epsilon)\frac12|xy|(1+\delta_\epsilon)$.
So $|xy|\leq |xn|+|yn| \leq |xy|(1-\delta_\epsilon)(1+\delta_\epsilon)$, a contradiction.
So we have $|m_{xy},gm_{xy}| < \frac12\epsilon|xy|(1+\delta_\epsilon)$, showing that $m_{xy} \in X_{\cont}$.

Lastly, if the action is linear it is clear that $X_{\cont}$ is closed under addition and scalar multiplication.
\end{proof}

So the continuous part always forms a uniformly convex subaction.
Of course, it might be empty, or in the linear case, consist of just the origin.
However in some cases we can prove that the continuous part is non-trivial.
To do this we need the following lemma.
For open neighbourhoods of the identity $U,V\subseteq G$ we write $U<V$ if there is an open $W$ such that $WU\subseteq V$.
Every topological group has a large chain of open neighbourhoods.
\begin{lemma}\label{lem-large-chain-of-open-neighbourhoods}
Let $G$ be a topological group and $U\subseteq G$ an open neighbourhood of the identity.
Then there are open neighbourhoods of the identity $V_t$ for $t\in(0,1]$ with $V_1=U$ and $V_t<V_{t'}$ if $t<t'$.
\end{lemma}
\begin{proof}
Let $U_0 = U$ and for $n\geq 1$ we pick open neighbourhoods of the identity $U_n$ recursively with $U_n^2\subseteq U_{n-1}$.
For $t\in [0,1]$ let $\sigma_t$ be the set of finite sequences of integers $\{ (n_1,n_2,\ldots,n_s) \mid n_1>n_2> \cdots > n_s \geq 0 \text{ and } \sum_{i=1}^s2^{-n_i}\leq t\}$.
Define
\[V_t = \bigcup_{(n_1,n_2,\ldots,n_s)\in \sigma_t}U_{n_1}U_{n_2}\cdots U_{n_s}.\]
Note that when $t=2^{-n_1} + 2^{-n_2} +\cdots + 2^{-n_s}$ is a dyadic rational with $n_1>n_2>\cdots >n_s$ we have $V_t = U_{n_1}U_{n_2}\cdots U_{n_s}$.
In particular, $V_1=U$.

For $t\leq t'$ we have $V_t\subseteq V_{t'}$.
Moreover, there is a dyadic rational $2^{-n_1}+\cdots + 2^{-n_s}$ with $n_1>\cdots>n_s$ and $t<2^{-n_1}+\cdots + 2^{-n_s}<t'$.
There is also $n_0>n_1$ with $2^{-n_0}+2^{-n_1}+\cdots + 2^{-n_s} <t'$.
Then $U_{n_0}V_t \subseteq U_{n_0}U_{n_1}\cdots U_{n_s} \subseteq V_{t'}$.
So $V_t < V_{t'}$.
\end{proof}

Now we can show Proposition \ref{prop:continuousSubaction}, which for convenience, we restate below.
\begin{proposition}
Let $G$ be a topological group and let $U$ be an open neighbourhood of the identity. Let $X$ be a uniformly convex metric space.
Let $\alpha\colon G\curvearrowright X$ be an isometric action.
Let $c>0$ and suppose there is $x\in X$ such that $d(\alpha(g)x,x)\leq c$ for all $g\in U$.
Then there is $y\in X$ with $d(x,y)\leq 2c$, such that the restriction of $\alpha$ to the closed convex hull of $\alpha(G)y$ is continuous.
\end{proposition}

\begin{proof}
Let $(V_t)$ be as in Lemma \ref{lem-large-chain-of-open-neighbourhoods}.
We will consider the centres of the sets $\alpha(V_t)(x)$.
For $t\in(0,1)$ define $r(t) = \rho(\alpha(V_t)(x))$.
This is an increasing function, so there is some $t_0\in(0,1)$ where $r$ is continuous.
Define $y = Z(\alpha(V_{t_0})(x))$.
Since $|ax|\leq c$ for all $a\in\alpha(V_{t_0})(x)$, we have $r(\alpha(V_{t_0})(x))\leq c$ and $|xy|\leq 2c$.
We will show that the map $g\to \alpha(g)y$ is continuous.
Since $|\alpha(g)y,\alpha(h)y| = |\alpha(g^{-1}h)y,y|$, it is enough to show that $|\alpha(g)y,y|$ is small for $g$ close to the identity.

Let $\epsilon>0$ and $t>t_0$.
Let $\kappa_\epsilon$ be as in Lemma \ref{lem-additional-information-centre}(ii).
There is an open $W$ with $WV_{t_0}\subseteq V_t$.
Let $g\in W$.
We get
\[\alpha(g)y = Z(\alpha(g)\alpha(V_{t_0})(x)) = Z(\alpha(gV_{t_0})(x)).\]
We have $gV_{t_0} \subseteq V_t$.
Since $\rho(\alpha(g)\alpha(V_{t_0})(x)) = \rho(\alpha(V_{t_0})(x)) = r(t_0)$, we get
\[|\alpha(g)y , Z(\alpha(V_t)(x))|  \leq \epsilon r(t) + \kappa_\epsilon (r(t)-r(t_0))\]
and also
\[|y , Z(\alpha(V_t)(x))|\leq \epsilon r(t) + \kappa_\epsilon(r(t)-r(t_0))\]
so
\[|\alpha(g)y-y|\leq 2\epsilon r(t) + 2\kappa_\epsilon(r(t)-r(t_0)).\]
Since $r$ is continuous at $t_0$, the right-hand side goes to zero as $\epsilon\to 0$ and $t\to t_0$.
This shows that $g \to \alpha(g)y$ is continuous, so $y\in X_{\cont}$.
By Lemma \ref{lem:X-cont-subspace}, the restriction of $\alpha$ to the closed convex hull of $\alpha(G)y$ is continuous.
\end{proof}

Below is an alternative proof for locally compact groups acting on uniformly convex Banach spaces.
The existence of a Haar measure allows us to use a regularization procedure based on our notion of mean centre. 
\begin{proof}
Assume that $G$ is locally compact and that $X$ is a uniformly convex Banach space.
Let $K\subseteq U$ be a compact neighbourhood of the identity.
Let $\mu$ be the left-invariant Haar measure.
Let $f\colon U \to X$ be the orbit map $f(a) = \alpha(a)x$.
We will use it to push-forward the Haar measure.
Let $y = \tilde Z((f_{|K})_*\mu)$ be the mean centre in $X$ (see Construction \ref{constr:meancentre}).
Then $\norm{x-y}\leq 2c$.
We will show that $\norm{\alpha(g)y-y}$ is small for $g$ close to the identity.

Let $\epsilon>0$.
Since $\mu$ is outer regular, there is $K\subseteq U'\subseteq U$ with $\mu(U')\leq (1+\epsilon)\mu(K)$.
Let $V$ be an open neighbourhood of the identity with $VK \subseteq U'$.
Let $g\in V$.
Let $\lambda_g\colon K\to gK$ be left multiplication by $g$.
Then $\lambda_g$ is a measure-preserving map and $\alpha(g)\circ f_{|K} = f_{|gK} \circ \lambda_g$ so $\alpha(g)y=\tilde Z((\alpha(g)\circ f_{|K})_*\mu) = \tilde Z((f_{|gK})_*\mu)$.
We get
\begin{align*}
    \alpha(g)y-y &= \tilde Z((f_{|gK})_*\mu)-\tilde Z((f_{|K})_*\mu)\\
    &= \frac{\mu(gK\setminus K)}{\mu(K)}(\tilde Z((f_{|gK\setminus K})_*\mu) - \tilde Z((f_{|K\setminus gK})_*\mu)).
\end{align*}
So
\[\norm{\alpha(g)y-y}\leq \frac{4\mu(K\setminus gK)}{\mu(K)}\rho(\alpha(U)x)\leq \frac{4c(\mu(U')-\mu(K))}{\mu(K)}\leq 4c\epsilon.\]
This shows that the map $g\to \alpha(g)y$ is continuous near the identity, so $y\in X_{\cont}$.
By Lemma \ref{lem:X-cont-subspace}, the restriction of $\alpha$ to the closed convex hull of $\alpha(G)y$ is continuous.
\end{proof}
This theorem gives us some corollaries such as the ones mentioned in the introduction. To state more of them, we need first to introduce good notions of isometric actions for bornological groups.

\subsection{Isometric actions of bornological groups}\label{sec:FX}

Given a bornological group $G$, a metric space $X$, and an action $\sigma$ of $G$ by isometries on $X$, we shall consider two natural notions of ``compatibility'' of the action with the bornology on $G$.
\begin{definition}
Let $G$ be a bornological group and let $\sigma$ be an action of $G$ by isometries on a metric space $X$. We say that the action is controlled (resp.\ well controlled) if the orbits of bounded sets are bounded (resp.\ are relatively compact).
\end{definition}
Note that a unitary representation is trivially controlled, but not necessarily well controlled.
Given a bornological group $G$ and a metric topological space $X$, a map $\varphi:G\to X$ is called proper if the preimage of a bounded set is bounded.

\begin{definition}
Let $G$ be a bornological group and let $\sigma$ be an action of $G$ by isometries on a metric space $X$. We that $\sigma$ has almost fixed points if for every $\eps>0$ and every bounded subset $B\subset G$, there exists $x\in X$ such that $d(\sigma(g)x,x)\leq \eps$ for all $g\in B$.
\end{definition}
\begin{definition}\label{def:T_B}
Let $\pi$ be a norm-preserving representation of a bornological goup $G$ on a Banach space $E$.  We that $\pi$ 
 has almost invariant vectors if the isometric action induced by $\pi$ on the unit sphere has almost fixed points.
\end{definition}

\begin{definition}
Let $G$ be a bornological group, let $\mathcal X$ be a family of metric spaces, and let $\mathcal E$ be a family of Banach spaces.
\begin{enumerate}[(i)]
\item We say that $G$ has coarse property F$\mathcal{X}$ (resp.\ F$\mathcal X^-$) if every controlled (resp.\ well-controlled) isometric action of $G$ on a space from $\mathcal X$ has a fixed point. 

\item We say that $G$ has coarse property (T$\mathcal E$) (resp.\ (T$\mathcal E^-$)) if all its norm-preserving representations (resp.\ well-controlled representations) on elements of $\mathcal E$ with almost invariant vectors have a non-zero invariant vector.
\end{enumerate}
\end{definition}

We are now ready to list some consequences of Proposition \ref{prop:continuousSubaction}, which in particular applies to  locally compact groups. 

\begin{corollary}\label{cor:THetc}
Let $\mathcal E$ be a family of uniformly convex Banach spaces closed under taking closed linear subspaces.
Let $\mathcal X$ be a family of uniformly convex metric spaces closed under taking convex closed subsets.
Let $(G,\tau,\mathcal B)$ be a locally bounded topological bornological group.
\begin{enumerate}[(i)]
    \item The topological bornological group $(G,\tau,\mathcal B)$ has topological property (T$\mathcal E$) if and only if the bornological group $(G,\mathcal B)$ has coarse property (T$\mathcal E$).
    \item The topological bornological group $(G,\tau,\mathcal B)$ has the topological fix-point property F$\mathcal{X}$ if and only if the bornological group $(G,\mathcal B)$ has the coarse fix-point property F$\mathcal{X}$.
\end{enumerate}
\end{corollary}
\begin{proof}
\begin{enumerate}[(i)]
\item Suppose $(G,\tau,\mathcal B)$ has topological property (T$\mathcal E$).
Let $E\in\mathcal E$.
Consider a representation $\pi\colon G\curvearrowright E$ with almost invariant vectors.
Let $U$ be an open bounded set in $G$.
Let $A\in\mathcal B$ and $\epsilon>0$ and let $v\in E$ be a unit vector with $\norm{gv-v}\leq\epsilon$ for all $g\in U\cup A$.
By Proposition \ref{prop:continuousSubaction}, there is $w\in E_{\cont}$ with $\norm{v-w}\leq2\epsilon$.
In particular $\norm w\geq1-2\epsilon$, and for $g\in K$ we have $\norm{gw-w}\leq\norm{gw-gv}+\norm{gv-v}+\norm{v-w}\leq5\epsilon$, so $\norm{g\frac w{\norm w}-\frac w{\norm w}}\leq\frac{5\epsilon}{1-2\epsilon}$.
This shows that $E_{\cont}$ has almost invariant vectors, so it has an invariant vector, showing that $(G,\mathcal B)$ has coarse property $(T\mathcal E)$.
The other direction is trivial.

\item Suppose $(G,\tau,\mathcal B)$ has topological property F$\mathcal{X}$.
Let $\alpha\colon G\curvearrowright X$ be a controlled isometric action on an element of $\mathcal X$.
Since $\alpha$ is controlled, we can apply Proposition \ref{prop:continuousSubaction} on any element of $X$ to find that $X_{\cont}$ is non-empty.
Then $X_{\cont}$ has a fixed point, so $(G,\mathcal B)$ has coarse property F$\mathcal{X}$.
The other direction is trivial.
\end{enumerate}
\end{proof}
\begin{corollary}\label{cor:Tdense}
Let $\mathcal E$ be a family of uniformly convex Banach spaces closed under taking closed linear subspaces.
Let $\mathcal X$ be a family of uniformly convex metric spaces closed under taking convex closed subsets.
Let $G$ be a compactly generated locally compact group, and $H$ be a dense subgroup, and let $U$ be a relatively compact open symmetric generating subset $G$. \begin{enumerate}[(i)]
    \item Assume $G$ has property $(T\mathcal E)$. Then any norm-preserving representation of $H$ on a space $E\in \mathcal E$ with $(U\cap H)$-almost invariant vectors has an invariant vector.
    \item Assume $G$ has property F$\mathcal{X}$. Then any affine isometric action of $H$ on some space $X\in \mathcal X$ such that $(U\cap H)$-orbits are bounded has a fixed point.
\end{enumerate}
\end{corollary}
\begin{proof} We equip $H$ with the topology $\tau_G$ induced by the topology on $G$, and the coarse structure $\mathcal K_G$ induced by the coarse structure $\mathcal K$ on $G$. 
This makes $(H,\tau_G,\mathcal K_G)$ is a locally bounded topological bornological group. Since $U$ is open, and $H$ is dense, for all $n\geq 1$, $(U^n\cap H)=(U\cap H)^n$. But since $U$ is symmetric and generates $G$, for any compact $K\subset G$, there exists $n\in \mathbb N$ such that $K\subset U^n$. Combining these two remarks, we get that $K\cap H\subset (U^n\cap H)$. Hence $U^n\cap H$ generates the bornology of $(H,\mathcal K_G)$. Moreover, since $H$ is dense in $G$, any continuous isometric action of $H$ on some complete metric space extends to $G$.  Hence topological $(T\mathcal E)$ and topological $F\mathcal X$ for $(H,\tau_G,\mathcal K_G)$ are equivalent respectively to $(T\mathcal E)$ and $F\mathcal X$ for $G$. Now we deduce from Corollary \ref{cor:THetc} that $(H,\tau_G,\mathcal K_G)$ has coarse $(T\mathcal E)$ under the assumption of (i), and  that $(H,\tau_G,\mathcal K_G)$ has coarse $F\mathcal X$ under the assumption of (ii). We conclude from the fact that the assumption of (i) that the representation of the bornological group $(H,\mathcal K_G)$ has almost invariant vectors, while the assumption of (ii) implies that the action  of $(H,\mathcal K_G)$ is controlled.
\end{proof}

 \section{Applications of Proposition \ref{prop:continuousSubaction} for locally compact groups}

The main goal of this section is to prove applications of Proposition \ref{prop:continuousSubaction} such as Corollary \ref{cor:Tcompactlypresented}, as well as others in the same spirit. In order to do so, we need to consider ultraproducts of isometric actions, which will be better defined in the setting of bornological groups. 

\subsection{Ultraproduct of isometric actions of bornological groups}

Let $(G,S)$ be a group equipped with a symmetric generating subset $S$. 
A controlled isometric affine action of $(G,S)$ on a metric space $X$ is an action whose $S$-orbits are bounded.
Now, let $G_i$ be an $I$-indexed net of groups and for each $i\in I$, let $S_i$ be a generating set of $G_i$.
Let $\mathcal{U}$ be a free ultrafilter on $I$, and denote $(G_\mathcal{U},S_\mathcal{U})$ the bounded ultraproduct of the net $(G_i,S_i)$.
 The displacement of a controlled isometric action $\sigma$ of $(G,S)$ is the infimum over all $x\in X$ of $d(x):=\sup_{s\in S}d(\sigma(s)x,x)$. 
We have
 \[d(x)-d(y)=\sup_{s\in S}d_X(\sigma(s)x,x)-\sup_{s\in S}d_X(\sigma(s)y,y)\leq \sup_{s\in S}(d_X(\sigma(s)x,x)-d_X(\sigma(s)y,y)).\]
 By triangular inequality, we have
 \[d_X(\sigma(s)x,x)\leq d_X(\sigma(s)x,\sigma(s)y)+d_X(\sigma(s)y,y)+d_X(y,x)=2d_X(y,x)+d_X(\sigma(s)y,y).\]
 We deduce that
 \[d(x)-d(y)\leq 2d_X(x,y).\]
 Hence $d:X\to \mathbb{R}_+$ is $2$-Lipschitz, hence continuous.
\begin{definition}
We let $\mathcal X$ be a family of metric spaces. We say that $\mathcal X$ is stable under rescaled ultraproducts if for every $I$-indexed net of triples $(X_i,\lambda_i,x_i)$, where $x_i\in X_i$ and $\lambda_i>0$, every ultrafilter on $I$, the ultraproduct of pointed metric spaces $(X_i,\lambda_id_{X_i},x_i)$ belongs to $\mathcal X.$
\end{definition}

\begin{lemma}\label{lem:displacement}
Let $\sigma$ be an isometric action without fixed points of $(G,S)$ on $X$. Then for all $\alpha<1$ and all $R\geq 0$, there exists $x$ such that $\inf_{|xy|\leq Rd(x)}d(y)\geq \alpha d(x) $.
\end{lemma}
\begin{proof}
Assume by contradiction that there exists $\alpha<1$ and $R\geq 0$ such that for all $x$, \[\inf_{|xy|\leq Rd(x)}d(y)< \alpha d(x).\] 
Fix some $x_0$ and choose recursively $x_{n+1}\in B(x_n,Rd(x_n))$ such that $d(x_{n+1})\leq \alpha d(x_{n}) $. By an immediate induction we get $d(x_n)\leq \alpha^n d(x_0)$. This implies that $|x_{n+1}x_n|\leq R\alpha^n$. Therefore $(x_n)$ is a Cauchy sequence, and let $x$ be its limit. We deduce by continuity of $d$ that $d(x)=0$, meaning that $x$ is a fixed point: contradiction.
\end{proof}

Given some isometric actions without fixed points of groups $G_i$, we will construct an isometric action with positive displacement (in particular without fixed point) on an ultraproduct of these groups.
Since we do not want to restrict to separable metric spaces, we will have to consider nets in full generality, rather than just sequences.
For technical reasons we can only use ultrafilters that are not $\sigma$-complete.
Recall that an ultrafilter $\mathcal U$ on a set $I$ is not $\sigma$-complete precisely when there is a function $f\colon I\to \mathbb N$ such that $f^{-1}(n) \not\in\mathcal U$ for all $n\in \mathbb N$.
If there are no measurable cardinals, then every non-principal ultrafilter is not $\sigma$-complete, and this is consistent with ZFC.
Actually, we only need the existence of suitable not $\sigma$-complete ultrafilters and we can just prove this.

\begin{lemma}\label{lem-ultrafilter_on_net_converging_to_infty}
Let $(I,\leq)$ be a directed set without maximum.
There is a net $n_i$ indexed on $I$ with values on $\N$, and a cofinal ultrafilter $\mathcal U$ such that $\lim_\mathcal Un_i = \infty$.
\end{lemma}
\begin{proof}
We will first prove there are $n_i$ such that for all $N$ and all $i$ there is some $j\geq i$ with $n_j\geq N$.
For each $i\in I$, we choose an element $t(i) \in I$ which is larger.
This gives $I$ the structure of a directed graph, with an arrow from $i$ to $t(i)$.
Choose a set $A\subseteq I$ that contains exactly one element of each connected component of the graph.
For $i\in t^n(A)$, we put $n_i = n$.
Since $A$ contains only one element of each component, we are not trying to label an element twice, and it is well-defined.
Put $n_i = 0$ for all $i$ that did not get a number yet.
Now for any $i\in I$ there is $a\in A$ in the same component.
Then there is a path from $a$ to $i$ in the graph (possibly taking edges in the opposite direction), and this path contains a vertex that is above $a$ and $i$.
Then for any large enough $N$, we see that $t^N(a) \geq i$, and $n_{t^N(a)} = N$.

Consider the sets $\{j\mid j\geq i\} \cap \{j\mid n_j\geq N\}$ for all $i\in I, N\in \N$.
By construction of the $n_i$ these sets are all non-empty, and they generate a filter $\mathcal F$.
Let $\mathcal U$ be any ultrafilter containing $\mathcal F$, then $\mathcal U$ is a cofinal ultrafilter converging to $\infty$.
\end{proof}

\begin{proposition}\label{prop:ultralimitAffine}
Let $\mathcal X$ be a family of metric spaces which is stable under rescaled ultraproducts.
Let $(G_i,S_i)_{i\in I}$ be a net of groups equipped with symmetric generating subsets.
For every $i\in I$, let $\sigma_i$ be a controlled isometric action without fixed points of $(G_i,S_i)$ on $X_i$. For each $i$, there exist a cofinal ultrafilter $\mathcal U$ on $I$, a space $X_\mathcal{U}\in \mathcal X$, and a controlled isometric action $\sigma_\mathcal{U}$ of the ultraproduct $(G_\mathcal{U},S_\mathcal{U})$ on $X_\mathcal{U}$ with positive displacement.
\end{proposition}

\begin{proof}
Let $n_i$ and $\mathcal U$ be as in Lemma \ref{lem-ultrafilter_on_net_converging_to_infty}.
By Lemma \ref{lem:displacement}, there exist points $x_i\in X_i$ such that \[\inf_{|x_iy|\leq n_id_i(x_i)}d_i(y)\geq (1-1/n_i) d_i(x_i)\]
We can consider the ultralimit $\sigma_\mathcal{U}=\lim_\mathcal{U}\sigma_n$ which defines a controlled isometric action $\sigma_\mathcal U$ of $G_\mathcal{U}$ on the metric space $X_\mathcal{U}=\lim_\mathcal{U}(X_i,d_{X_i}/d_i(x_i),x_i)$.
For any $y = (y_i) \in X_\mathcal U$ there is $k\in\N$ with $|x_iy_i| \leq kd_i(x_i)$ for all $i$.
For all $i$ with $n_i \geq k$ we get $d_i(y_i) \geq (1-1/n_i)d_i(x_i)$, so there is $s_i \in S_i$ with $|s_iy_i,y_i|\geq (1-2/n_i)d_i(x_i)$.
Now $(s_i) \in S_\mathcal U$ and $d_{X_\mathcal{U}}((s_i)\cdot (y_i),(y_i))  \geq 1$.
So $d(y) \geq 1$.
In particular $\sigma_\mathcal U$ has positive displacement.
\end{proof}

We can reformulate the proposition as follows.

\begin{corollary}\label{cor-ultraproduct-groups-without-FH}
Let $\mathcal X$ be a family of metric spaces which is stable under rescaled ultraproducts.
Let $G_i$ be a net of bornological groups such that the bornology of $G_i$ is generated by $S_i\subseteq G_i$.
Suppose that the $G_i$ do not have property (F$_\mathcal X$).
Then there is a cofinal ultrafilter $\mathcal U$ such that the ultraproduct $G_\mathcal U$, with bornology generated by $S_\mathcal U$, does not have property (F$_\mathcal X$). Better: $G_\mathcal U$ admits a controlled isometric action with positive displacement on some element of $\mathcal X$. 
\end{corollary}

See also \cite[Proposition 1.1]{CheCowStr04} for a similar result.

\subsection{Boundedly presented groups}
\begin{definition}\label{def-boundedly-presented}
A bornological group $G$ is \emph{boundedly presented} if there is a bounded symmetric set $S$ which generated the bornology, and a set of relations $R$ consisting of words in $S$ with a bounded length, such that $G = \langle S \mid R\rangle$.
\end{definition}

In this section, we show that every bornological group with FH is a quotient of a boundedly presented group with FH. More generally,
we make this statement for F$\mathcal X$ where $\mathcal X$ is any family of uniformly convex metric spaces closed under taking closed convex subspaces and rescaled ultraproducts.

\begin{theorem}\label{thm-boundedly-presented}
Let $\mathcal X$ be a family of uniformly convex metric spaces closed under taking closed convex subspaces and rescaled ultraproducts.
Let $G$ be a bornological group with bounded geometry and bornology generated by a single set.
Assume $G$ has property F$\mathcal{X}$.
There is a boundedly presented bornological group $G'$ with bounded geometry and property F$\mathcal{X}$ together with a quotient map $G' \twoheadrightarrow G$.
\end{theorem}
\begin{proof}
By assumption there is a symmetric subset $S\subseteq G$ generating the bornology, which we can also take to be gordo.
Then $S^2\subseteq g_1S \cup \cdots \cup g_k S$.
Let $N$ large enough that $g_i \in S^N$ for $1\leq i\leq k$.
Let $R$ be the set of all words in $S$ that evaluate to the identity in $G$, so that $G = \langle S\mid R\rangle$.
For $n \geq N + 3$, let $R_n$ be all words in $R$ of length at most $n$, and let $G_n = \langle S \mid R_n\rangle$.
Note that $G$ is a quotient of $G_n$.
Since $n\geq N+3$, we already have $S^2 \subseteq g_1S \cup \cdots \cup g_kS$ in $G_n$, showing that $G_n$ has bounded geometry.
Suppose that none of the $G_n$ have property F$\mathcal X$.
Then by Corollary \ref{cor-ultraproduct-groups-without-FH}, there is a free ultrafilter $\mathcal U$ on $\mathbb N$ such that the ultraproduct $G_\mathcal U$ does not have property F$\mathcal X$.
We have a map $f\colon G \to G_\mathcal U$ as follows: for each $g\in G$, write $g = s_1\cdots s_l$ with $l$ minimal, and put $f(g) = (s_1\cdots s_l) \in G_\mathcal U$.
The choice of representation of $g$ is not unique, but the element $s_1\cdots s_l \in G_n$ will be unique for large enough $n$, thus giving a well-defined map.
Now let $(g_n)\in G_\mathcal U$ be any element.
Note that the $g_n$ are all in $S^l$ for some fixed $l$, and $S^l$ is covered by finitely many left translates of $S$.
So there are $h_1,\ldots,h_m \in G$ and for each $n$ there is $1\leq i_n\leq m$ and $s_n\in S$ such that the equality $g_n = h_{i_n}s_n$ holds in $G$.
For large enough $n$, this inequality then holds in $G_n$ as well.
Then there is one $i$ such that $\{n\mid i_n = i\} \in \mathcal U$, and $(g_n) = f(h_i)\cdot (s_n)$.
So $G_\mathcal U = f(G)S_\mathcal U$, showing that $G$ is coarsely dense in $G_\mathcal U$.
By Proposition \ref{prop-dense-subgroup}, we get a contradiction.
So one of the $G_n$ satisfies our criteria.
\end{proof}

\subsection{Marked locally compact groups}
Since there exists no good notion of free group in the category of compactly generated locally compact groups, we shall rather work with a \emph{relative} notion of marked groups in this setting.

Given a locally compact compactly generated group, a generating subset $S$ will be called admissible if it is symmetric, compact, equal to the closure of its interior, and a neighbourhood of the identity.
\begin{definition}
Let $G^0$ be a locally compact compactly generated group and let $S^0$ be an admissible 
generating subset.
\begin{enumerate}[(i)]
    \item A $(G^0,S^0)$-marked group is a pair $(G,\pi_G)$, where $\pi_G:G^0\to G$ is a continuous surjective homomorphism whose kernel does not intersect the interior of $S^0$.
    \item An isomorphism $(G,\pi_G)\to (G',\pi_{G'})$ is a continuous isomorphism $\varphi:G\to G'$ such that  $\varphi\circ\pi_G = \pi_{G'}$. An isomorphism class is completely determined by the kernel of $\pi_G$.
    \item We denote by $\mathcal{M}(G^0,S^0)$ the set of isomorphism classes of $(G^0,S^0)$-marked groups.
  
     \item We equip $\mathcal{M}(G^0,S^0)$ with the Chabauty topology on closed subgroups of $G^0$.
\end{enumerate}
\end{definition}

Note that the subgroups $\ker_G$ are (uniformly) discrete, in particular they are closed. 
Let $H$ a closed subgroup of $G$. We recall a basis $(\Omega_{V,K}(H))_{K,V}$ of neighbourhoods of $H$ for the Chabauty topology, where  $V$ are neighbourhoods of $1_{G^0}$ and $K$ are compact subsets:
each set $\Omega_{V,K}(H)$  consists of all closed subgroups $H'$ such that $H\cap K\subset (H'\cap K)V$ and $H'\cap K\subset (H\cap K)V$.  
The space of closed subgroups of a locally compact group is well-known to be compact for the Chabauty topology.
Since the condition that the subgroups avoid the interior of $S^0$ is a closed condition, this turns $\mathcal{M}(G^0,S^0)$ into a compact topological space.

\begin{proposition}\cite[Proposition 8.A.15. L]{CdlH}
We let $G$ be a locally compact group, and let $S\subset G$ be an admissible generating set. Let $\hat{G}$ be the group defined by the presentation $\langle S\mid R\rangle $, where $R$ consists of all words of the form $s_1s_2s_3$ such that $s_1s_2s_3\equiv 1$ in $G$. Then $\hat{G}$ admits a unique locally compact group topology  such that the projection $\pi:\hat{G}\to G$ is a homeomorphism in restriction to $S$. \end{proposition}\label{prop:CornHarp}
\qed
\begin{corollary}\label{cor:denseCP}
We let $G$ be a locally compact group, and let $S\subset G$ be an admissible generating set. Then $(G,\pi_G)\in \mathcal{M}(G^0,S^0)$ for some compactly presented $G^0$. 
Moreover, every open set in $\mathcal M(G^0,S^0)$ containing $(G,\pi_G)$ contains an element $(G',\pi_{G'})$ such that $\ker(\pi_{G'})\subseteq \ker(\pi_G)$ and $G'$ is compactly presented.
 \end{corollary}
\begin{proof}
The first statement is an immediate consequence of Proposition \ref{prop:CornHarp}. For the second statement: let $H_n$ be the normal subgroup of $G^0$ spanned by $\ker \pi_G\cap (S^0)^n$ and let $G_n:=G^0/H_n$ for all $n\geq 2$. Clearly $H_n$ converges to $\ker \pi_G$ in the Chabauty topology.
By construction, we have $H_n=\ker \pi_{G_n}$.
This shows that $(G_n,\pi_{G_n})$  converges to $(G,\pi_G)$.
So we can put $(G',\pi_G' ) = (G_n,\pi_{G_n})$ for some large enough $n$.
Besides, a compact presentation for $(G',\pi_{G'})$ is given by $\langle S^0\mid R^0\cup R_n \rangle$, where $\langle S^0\mid R^0\rangle$ is a compact presentation of $G^0$ and $R_n$ is the set of words of length $\leq n$ in $S^0$ representing elements of $\ker\pi_G$.
\end{proof}

\begin{lemma}\label{lem:UltralimitToG}
Let $(G_i,\pi_{G_i})_{i\in I}\in \mathcal{M}(G^0,S^0)$ be a net converging to $(G,\pi)$, and let $\mathcal{U}$ be a cofinal ultrafilter on $I$. Denote $S_i=\pi_{G_i}(S)$, and let $(G_\mathcal{U},S_\mathcal{U})$ be the bounded ultraproduct with respect to $\mathcal{U}$ of $(G_i,S_i)$. There is a surjective homomorphism $p:G_\mathcal{U}\to G$ such that $p(S_\mathcal{U})=\pi_G(S)$. Moreover $\ker(p)\subset S_\mathcal{U}$.
\end{lemma}
\begin{proof}
Let $(G^0_\mathcal{U},S^0_\mathcal{U})$ denote the bounded ultraproduct of the constant sequence $(G^0,S^0)$ with respect to $\mathcal{U}$. Consider the canonical surjective homomorphism $\pi_\mathcal{U}=\lim_\mathcal{U}\pi_{G_i}:(G^0_\mathcal{U},S^0_\mathcal{U})\to (G_\mathcal{U},S_\mathcal{U})$.
Let $N^0$ be the subset of $G^0_\mathcal{U}$ consisting of sequences $g_i\in S^0$ such that $\lim_\mathcal{U}g_i=1_{G^0}$.

We note that $N=\pi_{\mathcal{U}}(N^0)$ coincides with the subset of $G_\mathcal{U}$ consisting of nets $\pi_{G_i}(g_i)\in S_i$ such that $\lim_\mathcal{U}g_i=1_G$. 
Since $G^0$ is locally compact, for all bounded nets $(g_i)$ in $G^0$, there exists a unique $g\in G^0$ such that $\lim_\mathcal{U}g_i=g$, or equivalently such that $\lim_\mathcal{U}g^{-1}g_i=\lim_\mathcal{U}g_ig^{-1}=1_{G^0}$. 
This implies that $N^0$ is a normal subgroup and that $G^0_\mathcal{U}/N^0$ is canonically isomorphic to $G^0$. It follows that $N$ is a normal subgroup of $G_\mathcal{U}$.
We have the following commutative diagram, which we want to complete by adding the arrow $p$.
\[
\begin{CD}
N^0                      @>>>  G_\mathcal{U}^0 @>p^0>> G^0         \\
@VVV               @V\pi_\mathcal{U} VV  @VV\pi_G V\\
N    @>>>    G_\mathcal{U} @>p?>> G
\end{CD}
\]

Now for all bounded nets $(k_i)\in \Pi_iG_i$, let $(g_i)$ be a bounded net in $G^0$ such that $\pi_{G_i}(g_i)=k_i$, and let $k=\pi_G(g)$ where $g=\lim_\mathcal{U}g_i$. We need to show that $k$ does not depend on the choice of lifted net $(g_i)$, or equivalently, that if $(h_i)$ is a bounded net such that $h_i\in \ker \pi_{G_i}$, then $\pi_G\circ p^0((h_i))=1_G$. By assumption, $\ker \pi_{G_i}$ converges to $\ker\pi_G$, hence $\lim_\mathcal{U}h_i\in \ker \pi_{G}$, so this is proved. 
Setting $p((g_i))=k$, this shows that $p$ is a well-defined morphism. 

It is clear by construction that $p$ is surjective (since $p^0$ is surjective). 
We are left to checking that $N=\ker p$. The inclusion $N\subset \ker p$ follows by commutativity of the diagram. For the other inclusion, let $(k_i)\in \ker p $. Let $(g_i)$ be a lift of $(k_i)$ in $G_\mathcal{U}^0$. We have that $p^0((g_i))$ belongs to the kernel of $\pi_G$, meaning that $\lim_\mathcal{U}g_i\in \ker \pi_G$. 
But since $\ker \pi_{G_i}$ converges to $\ker\pi_G$, this means that there exists $h_i\in \ker \pi_{G_i}$ such that $\lim_\mathcal{U}h_i^{-1}g_i=1_{G^0}$. On replacing $(g_i)$ by $(h_i^{-1}g_i)$, which is another lift of $(k_i)$, we have that $(g_i)\in N^0$. Hence $(k_i)$ lies in $\pi_\mathcal{U}(N^0)=N$. This ends the proof of the fact that $N=\ker p$.
 
Since $N\subset S_\mathcal{U}$, we also get the last statement of the lemma. 
\end{proof}

We recall for convenience the statement of Theorem \ref{thm-positive-displacement} before proving it.
\begin{theorem}\label{thmInSection:positive-displacement}
Let $\mathcal X$ be a family of uniformly convex metric spaces, closed under taking rescaled ultraproducts and closed under taking closed convex subspaces.
Let $G$ be a compactly generated locally compact group such that every continuous isometric action on an element of $\mathcal X$ as almost fixed points. Then $G$ has property F$\mathcal{X}$.
\end{theorem}
\begin{proof}
Let us argue by contradiction and assume $G$ does not have property F$\mathcal X$.  By Corollary  \ref{cor-ultraproduct-groups-without-FH}, $G_{\mathcal U}$ has a controlled action $\sigma$ on some $X\in \mathcal X$ with positive displacement. We let $p:G_\mathcal{U}\to G$ be the projection from Lemma \ref{lem:UltralimitToG}. Since $\ker(p)\subset S_\mathcal{U}$ and the action is controlled, we deduce that the orbits of  $\ker p$ are bounded. Hence the closed convex set $Y$ of $\ker p$-fixed point is non-empty and $\sigma(G_\mathcal{U})$-invariant. Moreover the restriction of $\sigma$ to $Y$ factors through $G$. Applying Proposition \ref{prop:continuousSubaction}, we find a closed convex subspace $Z\subset Y$ on which the action of $G$ is continuous. Since this action is a subaction of $\sigma$ it still has positive displacement: contradiction.
\end{proof}

Another corollary of Lemma \ref{lem:UltralimitToG} is the following. 

\begin{theorem}\label{thm:Topen}
Let $\mathcal X$ be a family of uniformly convex metric spaces which is stable under rescaled ultraproduct and closed under taking closed geodesic subspaces. 
The subset of $S$-marked groups  $(G,\pi_G)$ such that $G$ has property F$\mathcal{X}$ 
forms an open subset of $\mathcal{M}(G^0,S^0)$. \end{theorem}

\begin{corollary}
The subset of $S$-marked groups  $(G,\pi_G)$ such that $G$ has property (T) 
forms an open subset of $\mathcal{M}(G^0,S^0)$. 
\end{corollary}
\begin{proof}
We shall prove that groups without F$\mathcal{X}$ form a closed subset of $\mathcal{M}(G^0,S^0)$.
Let $(G_i,\pi_{G_i})\in \mathcal{M}(G^0,S^0)$ be a net converging to $(G,\pi_G)$, and let $\mathcal{U}$ be a cofinal ultrafilter as in Proposition \ref{prop:ultralimitAffine}. Denote $S_i=\pi_{G_i}(S)$.
 We deduce from Corollary \ref{cor-ultraproduct-groups-without-FH}  that the bounded ultraproduct $(G_\mathcal{U},S_\mathcal{U})$ of the net $(G_i,S_i)$ does not have  F$\mathcal{X}$.
Let $\sigma$ be a controlled isometric action without fixed points of the bornological group $(G_\mathcal{U},S_\mathcal{U})$ on some $X\in \mathcal X$. 
  Since $S_\mathcal{U}$ has bounded orbits, the kernel $N$ of the surjective homomorphism $p:G_\mathcal{U}\to G$ given by Lemma \ref{lem:UltralimitToG} has bounded orbits, hence admits a fixed point $v$. Since $N$ is normal in $G_\mathcal{U}$, the set $X_N$ of 
 $\sigma(N)$-fixed points is  a closed   $\sigma(G_\mathcal{U})$-invariant subspace of $X$.
 We observe that in a uniformly convex metric space, any pair of points is joined by a unique geodesic, hence it follows that the set of $\sigma(N)$-fixed points is convex, and therefore $X_N\in \mathcal X$. The restriction of $\sigma$ to $X_N$ induces an action of $G$, and by Proposition \ref{prop:continuousSubaction}, we may find a closed convex subspace of $X_N$ on which the action of $G$ is continuous. Hence $G$ does not have property F$\mathcal{X}$ and we are done.
\end{proof}

Finally, we restate Theorem \ref{thm:FCAT(0)compactlypresented} before proving it.
\begin{theorem}\label{thmInSection:FCAT(0)compactlypresented}
Let $\mathcal X$ be a family of uniformly convex metric spaces, closed under taking rescaled ultraproducts and closed under taking closed convex subspaces.
Let $G$ be a compactly generated locally compact group with property F$\mathcal{X}$. Then there exists a compactly presented locally compact group $G'$ with property F$\mathcal{X}$, and a continuous surjective homomorphism $G'\to G$ with discrete kernel.
\end{theorem}
\begin{proof}
 Let $S\subset G$ be a compact symmetric generating neighbourhood of the identity.
We then conclude by applying Corollary \ref{cor:denseCP} and Theorem \ref{thm:Topen}.
\end{proof}

\section{A focus on actions on Hilbert spaces}
 
The main goal of this section is to relate our two versions of coarse property (T) and coarse property FH, and to prove Proposition \ref{prop:wellcontrolledSubaction}. The first one is an adaptation the classical argument by Delorme and Guichardet, and while Proposition \ref{prop:wellcontrolledSubaction} might hold for actions on more general Banach spaces, we were only able to prove it for Hilbert spaces.

\subsection{Delorme-Guichardet for bornological groups}

Just like in the topological case, there is a connection between properties (T) and FH.
The proof is similar to the proof of the Delorme-Guichardet Theorem, as seen in \cite{BHV07}.
We say the bornology of a bornological group $G$ is generated by a single set if there is some bounded set $A$ such that for each bounded set $C\subseteq G$, there is some $n$ with $C\subseteq A^n$.
In particular, $A$ generates $G$ as a group.

\begin{proposition}\label{prop-Delorme-Guichardet}
Let $G$ be a bornological group.
\begin{enumerate}[(i)]
    \item If $G$ has coarse property (T-) (resp. (T)) then it has coarse property FH- (resp. FH).
    \item If $G$ has coarse property FH- (resp. FH) and its bornology is generated by a single set, then it has coarse property (T-) (resp. (T)).
\end{enumerate}
\end{proposition}
\begin{proof}
First suppose that $G$ has coarse property (T-).
As in \cite[Proposition 2.11.1]{BHV07}, we ca  construct a real Hilbert space $\H_t$, a map $\Phi_t\colon \H\to \H_t$ and an orthogonal representation $\pi_t\colon G\to O(\H_t)$ satisfying
\[\langle \Phi_t(v),\Phi_t(w)\rangle = \exp(-t\norm{v-w}^2)\]
and
\[\pi_t(g)\Phi_t(v) = \Phi_t(\alpha(g)v)\]
and the linear span of $\Phi_t(\H)$ is dense in $\H_t$.

Consider the direct sum of all these representations $\bigoplus_{t>0}\pi_t\colon G\to O(\bigoplus_{t>0}\H_t)$.
This has almost invariant vectors: let $A\subseteq G$ be bounded and $\epsilon>0$.
Consider the vector $\Phi_t(0)$ for some $t>0$.
It satisfies
\[\langle \pi_t(a)\Phi_t(0),\Phi_t(0)\rangle =\langle\Phi(t)(\alpha(a)(0)),\Phi_t(0)\rangle = \exp(-t\norm{\alpha(a)(0)}^2)\]
for each $a\in A$.
Since $\alpha$ is controlled, we can choose $t>0$ such that the right-hand side is always at least $1-\epsilon$.
Then $\norm{\pi_t(a)\Phi_t(0)-\Phi_t(0)}^2 = 2-2\langle\pi_t(a)\Phi_t(0),\Phi_t(0)\rangle \leq 2\epsilon$, showing that $\Phi_t(0)$ is an almost invariant vector.

Since $G$ has coarse property (T-) there is some invariant vector.
This has to have a non-zero coordinate in some $\H_t$, so there is some $t>0$ for which a non-zero invariant vector $w\in \H_t$ exists.
Suppose the orbits of $\alpha$ are unbounded.
Let $v\in \H$ be a vector and let $g_n\in G$ be a sequence such that $\{\alpha(g_n)v\}$ is unbounded.
Then the sequence $(\Phi_t(\alpha(g_n)v))$ converges weakly to 0.
For each $n$ we have
\[\langle w, \Phi_t(v)\rangle = \langle w,\pi_t(g_n)(\Phi_t(v))\rangle = \langle w,\Phi_t(\alpha(g_n)v)\rangle.\]
Taking the limit as $n\to\infty$ shows $\langle w,\Phi_t(v)\rangle = 0$.
Since $\Phi_t(\H)$ is dense in $\H_t$ it follows that $w=0$.
This is a contradiction, so $\alpha$ has bounded orbits, and has a fixed point as well.
So $G$ has property FH.

Suppose $G$ has property (T-).
We use the same construction as above, but we start with an affine action $\alpha$ that is well controlled, and we need the representation $\pi_t$ to be controlled.
So let $A\subseteq G$ be bounded and $w\in \H_t$.
We need to show that $\pi_t(A)w$ is relatively compact.
First consider the case that $w = \Phi_t(v)$ for some $v\in \H$.
The map $\Phi_t\colon\H\to \H_t$ is continuous, and $\alpha(A)v$ is relatively compact, so $\pi(A)w = \Phi_t(\alpha(A)v)$ is relatively compact as well.
Now if $w$ is in the span of $\Phi_t(\H)$ it follows that $\pi(A)w$ is still relatively compact.
Finally, since the span of $\Phi_t(\H)$ is dense in $\H_t$ it follows that $\pi(A)w$ is relatively compact for all $w\in\H_t$.
The rest of the argument is the same as above and it follows that $G$ has coarse property FH-.

Now suppose the bornology of $G$ is generated by some bounded set $A$, and that $G$ has coarse property FH.
Let $\pi\colon G\to O(\H)$ be an orthogonal representation to a real Hilbert space $\H$, and suppose that it has almost invariant vectors, but no invariant vector.
Then for any natural number $n$ we can find a unit vector $v_n\in \H$ with $\norm{\pi(a)v_n-v_n}\leq 2^{-n}$ for all $a\in A$.
Now define the map $b\colon G\to \bigoplus_n\H$ by $b(g) = (\pi(g)v_n-v_n)$.
We have $b(gh)=\pi(g)b(h)+b(g)$, and for $a\in A$, we have $\norm{b(a)}<1$ by the inequality above.
If $C\subseteq G$ is bounded, there is some $n$ with $C\subseteq A^n$, and it follows that $\norm{b(c)}<n$ (in particular, the norm is finite, and $b$ is well-defined).

Now define the affine action $\alpha\colon G\curvearrowright \bigoplus_n\H$ by $\alpha(g) = \bigoplus_n\pi(g) + b(g)$.
We see that $\alpha$ is controlled, so it has a fixed point $(w_n)$.
For large enough $n$ we must have $\norm{w_n}<1$.
For such $n$ we have $\alpha(w_n) = w_n$, so $w_n = \pi(g)(w_n) + \pi(g)(v_n)-v_n$ for all $g\in G$.
This gives $\pi(g)(w_n+v_n) = w_n + v_n$.
But $\pi$ has no fixed point, so $w_n = -v_n$, which as norm 1, a contradiction.
So $G$ has coarse property (T).

Lastly, suppose the bornology of $G$ is generated by some bounded set $A$, and that $G$ has coarse property FH-.
We use the same construction again, and we just have to show that if $\pi$ is controlled, then $\alpha$ is well controlled.
So let $(\eta_n)\in \bigoplus_n\H$ and let $C\subseteq G$ be bounded.
For each $n$, the set $\{\pi(C)\eta_n\}$ is relatively compact.
So the set $\{(\pi(c)\eta_n)\mid c\in C\}$ is also relatively compact.
The set $\pi(C)\xi_n$ is also relatively compact for each $n$, so $\{(\pi(c)(\xi_n)-\xi(n)\mid c\in C\}$ is relatively compact as well.
It follows that $\alpha(C)(\eta_n)=\{(\pi(c)\eta_n + \pi(c)(\xi_n) - \xi(n)\mid c\in C\}$ is relatively compact.
So $\alpha$ is well controlled.
The rest of the argument is the same as before, and it follows that $G$ has coarse property (T-).
\end{proof}

For bornological groups with bounded geometry, coarse property (T-) implies already that the bornology is generated by a bounded set.

\begin{proposition}\label{prop-generated-bounded-set}
Let $G$ be a bornological group with bounded geometry.
Suppose $G$ has coarse property (T-).
Then the bornology is generated by a bounded set.
\end{proposition}
\begin{proof}
Let $(A,\epsilon)$ be a Kazhdan pair.
We can assume that $A$ is also gordo.
Let $H$ be the subgroup generated by $A$.
Consider the representation $\pi \colon G\to U(l^2(G/H))$, given on basis vectors $e_{gH}$ by $\pi(g')e_{gH} = e_{g'gH}$.
We show that this is a controlled representation.
Let $C\subseteq G$ be bounded and consider a basis vector $e_{gH} \in l^2(G/H)$.
since $A$ is gordo, we can find $g_1,\ldots,g_n$ with $Cg\subseteq g_1A\cup \cdots \cup g_nA$.
Then $\pi(C)e_{gH}  = \{e_{cgH} \mid c\in C\} \subseteq \{e_{g_iH} \mid 1\leq i\leq n\}$ is finite.
It follows that $\pi(C)\xi$ is compact for all $\xi \in l^2(G/H)$, so $\pi$ is a controlled representation.

The basis vector $e_H\in l^2(G/H)$ satisfies $\pi(a)e_H = e_H$ for all $a\in A$.
Since $(A,\epsilon)$ is a Kazhdan pair, it follows that the representation has a constant vector.
But this means that $G/H$ is finite.
So if we add finitely many elements to $A$, we can generate the whole group with a bounded set $A'$.

To finish the proof we need to use bounded geometry again.
For a bounded set $C$ we can write $C\subseteq g_1A\cup\cdots\cup g_nA$.
We can find $k$ such that $g_i \in A^k$ for $1\leq i\leq n$.
Then $C\subseteq A^{k+1}$.
So the bornology is generated by a bounded set.
\end{proof}

\subsection{Proof of Theorem \ref{thm-(T)-and-(T-)-equivalent}}

Let us recall the statement of Proposition \ref{prop:wellcontrolledSubaction}.

\begin{proposition}\label{propInSection:wellcontrolledSubaction}
Let $G$ be a bornological group with bounded geometry and let $\alpha \colon G\curvearrowright \H$ be a controlled affine action on a Hilbert space.
There is a non-empty invariant closed subspace $\H'$ such that the restriction of $\alpha$ to $\H'$ is well controlled.
\end{proposition}
\begin{proof}
Let $C\subseteq G$ be gordo.
Let $\Lambda = \alpha(C)(0)$.
Let $v$ be a shopping centre of $\Lambda$.
We will show that $\alpha(D)(v)$ is relatively compact for all bounded $D\subseteq G$.

Since $C$ is gordo we can write $DC\subseteq g_1C\cup g_2C\cup\cdots\cup g_kC$, with $g_i\in G$.
Let $\Delta = \alpha(DC)(0) \subseteq \bigcup_{i=1}^k\alpha(g_i)\Lambda$.
Every element of $\alpha(D)(v)$ can be written as $\alpha(g)v$ with $g\in D$, and this is a shopping centre of $\alpha(g)\Lambda\subseteq \Delta$.
By Lemma \ref{lem-mall-relatively-compact}, we see that $\alpha(D)(v)$ is relatively compact.

Now it follows easily that $\alpha(D)w$ is relatively compact for all $w\in \alpha(G)v$.
Then this also holds for $w$ in the affine hull of $\alpha(G)v$.
Let $\H'$ be the closure of the affine hull of $\alpha(G)v$.
Then it follows that the restriction of $\alpha$ to $\H'$ is well controlled.
\end{proof}
 \begin{theorem}\label{thmInSection-(T)-and-(T-)-equivalent}
For bornological groups with bounded geometry, (T) and (T-) are equivalent.
Similarly, FH and FH- are equivalent for these groups.
\end{theorem}
\begin{proof}
Suppose $G$ has FH-.
Consider a controlled affine action $\alpha\colon G\curvearrowright \H$.
By Proposition \ref{prop:wellcontrolledSubaction}, there is a non-empty closed invariant subspace $\H'$ such that the restriction of $\alpha$ to $\H'$ is well controlled.
Then there is a fixed point in $\H'$.
So $G$ has FH.

Now note that $G$ has property (T) if and only if its bornology is generated by a bounded set and $G$ has property FH, while $G$ has property (T-) if and only if its bornology is generated by a bounded set and $G$ has property FH-.
So (T) and (T-) are also equivalent.
\end{proof}

We will use this theorem to prove the following result. 
For a bornological group $G$ and an ultrafilter $\mathcal U$ on any index set $I$, we consider the bounded ultraproduct $G_\mathcal U$ which is the quotient of bounded sequences $\{(g_n) \mid \{g_n\}\text{ is bounded}\}$ by the equivalence relation $(g_n)\sim (h_n)$ if $\{n\mid g_n = h_n\} \in \mathcal U$.
This has a natural bornology generated by the images of $A^I$ for bounded $A\subseteq G$.
\begin{proposition}\label{prop-ultrafilters-dont-matter}
Let $G$ be a bornological group with bounded geometry and let $\mathcal U$ be an ultrafilter on $I$.
Then $G$ has property FH if and only if $G_\mathcal U$ has property FH.
\end{proposition}
\begin{proof}
Let $C\subseteq G$ be a gordo set.
Assume that $G$ has property (FH).
We have an embedding $j:G\to G_\mathcal U$.
Consider $(g_i)\in G_\mathcal U$.
Since all $g_i$ must be contained in a bounded set and $C$ is gordo, we can find finitely many $h_1,\ldots,h_n\in G$ such that for all $i$ there is $1\leq k_i\leq n$ and $c_i\in C$ with $g_i=h_{k_i}c_i$.
Then there is one $1\leq k\leq n$ such that $\{i\mid k_i = k\}\in\mathcal U$.
We find $(g_i) = h_k(c_i) \in j(G)C_\mathcal U$.
Thus $j(G)C_\mathcal U = G_\mathcal U$, showing that $G$ is coarsely dense in $G_\mathcal U$.
Now consider a controlled action $\alpha\colon G_\mathcal U\curvearrowright\H$ on a Hilbert space.
Precomposing with $j$ yields a controlled action of $G$, which must have a fixed point $v$.
Then the orbit $\alpha(G_\mathcal U)v$ is bounded, so its centre must be a fixed point for $\alpha$.
So $G_\mathcal U$ has property (FH).

Now suppose that $G_\mathcal U$ has property (FH).
We will show that $G$ has property (FH-), which will be enough by Theorem \ref{thm-(T)-and-(T-)-equivalent}.
So let $\alpha\colon G\curvearrowright\H$ be a well-controlled action.
We define the new representation $\beta\colon G_\mathcal U\curvearrowright\H$ by the formula
\[\beta((g_i))v = \lim_\mathcal U\alpha(g_i)v.\]
For any element $(g_i)$ of $G_\mathcal U$ the set $\{(g_i)\}$ is bounded, so the set $\{\alpha(g_i)v\}\subseteq \H$ is relatively compact, showing that the limit on the right hand side exists.
So $\beta$ is well-defined, and then it is straight-forward to check that it defines a well-controlled affine action.
Therefore $\beta$ has a fixed point $v$, which is then a fixed point for $\alpha$ as well.
\end{proof}

\subsection{Proof of Theorem \ref{thm:Serre}}
We need the following well-known lemma.

\begin{lemma}\label{lem:L}
Let $\pi\colon G\curvearrowright\H$ be a unitary representation of a group $G$. 
Denote by $p^G$ the projection on the $G$-invariant vectors.
Then $p^Gv$ belongs to the closed convex hull $\pi(G)v$ for all $v\in\H$.
\end{lemma}
\begin{proof}
There exists a unique vector $v_0$ of minimal norm in the closed convex hull of $\pi(G)v$, which is then $\pi(G)$-invariant.
For all $w\in \pi(G)v$ we have that $w-v$ is perpendicular to the $G$-invariant vectors.
So this holds for $v_0$ as well, showing that $v_0 = p^Gv$.
\end{proof}
\begin{theorem}\label{thmInSection:Serre}
Let $(G,\mathcal B)$ be a bornological group with bounded geometry, whose bornology is generated by a single set.
Let $p\colon G\to Q$ be a surjective morphism and equip $Q$ with the induced bornology $p(\mathcal B)$.
Assume that the kernel $A=\ker(p)$ is central in $G$ and coarsely contained in the commutator $[G,G]$.
Then $G$ has coarse property (T) if and only if $Q$ has coarse property (T).
\end{theorem}
\begin{proof}
If $G$ has property (T) then so does $Q$, since it is a quotient.
So suppose $Q$ has property (T), and therefore also property (FH).
Let $\mathcal U$ be a free ultrafilter on $\mathbb N$.
First we show that any affine action of the bounded ultraproduct $G_\mathcal U$ on a Hilbert space $\H$ has almost invariant vectors.
Note that we have a surjective morphism $G_\mathcal U\to Q_\mathcal U$ with central kernel $A_\mathcal U$.

So let $\sigma\colon G_\mathcal U\curvearrowright \H$ be an affine transformation.
By Proposition \ref{prop:wellcontrolledSubaction} we may assume that $\sigma$ is well controlled.
We write $\sigma(g) v = \pi(g)v + b(v)$ where $\pi$ is a representation and $b$ is a cocycle for $\pi$, i.e. it satisfies $b(gh) = \pi(g)b(h) + b(g)$ for $g,h\in G_\mathcal U$.

Let $\H^G$ be the subspace of $\H$ consisting of $\pi(G_\mathcal U)$-invariant vectors and let $b^G$ denote the projection of $b$ onto $\H^G$.
Then the cocycle condition gives $b^G(gh) = b^G(g) + b^G(h)$, so $b^G\colon G_\mathcal U\to \H^G$ is a morphism to an abelian group.
Since $A$ has bounded geometry, we can prove similarly to the proof of \ref{prop-ultrafilters-dont-matter} that $A$ is coarsely dense in $A_\mathcal U$.
Since $A$ is coarsely contained in $[G,G]$ it follows that $A_\mathcal U$ is also coarsely contained in $[G_\mathcal U,G_\mathcal U]$.
Since $b^G$ is a controlled map we find that $b^G$ is bounded on $A_\mathcal U$.
But then for any $a\in A_\mathcal U$ and $n\in\mathbb Z$ we have $b^G(a^n) = nb^G(a)$ and boundedness implies that $b^G$ is 0 on $A_\mathcal U$.

Now let $\H^A$ be the subspace of $\H$ consisting of $\pi(A_\mathcal U)$-invariant vectors and let $p^A$ be the projection onto $\H^A$ and $b^A = p^A\circ b$.
The projection of the cocycle condition gives $b^A(gh) = \pi(g)b^A(h) + b^A(g)$ for $g,h\in G_\mathcal U$.
Using that $A_\mathcal U$ is central in $G_\mathcal U$ we find, for $a\in A_\mathcal U$ and $g\in G_\mathcal U$:
\[b^A(g) + b^A(a) =  \pi(a)b^A(g) + b^A(a) = b^A(ag) = b^A(ga) = \pi(g)b^A(a) + b^A(g)\]
so $b^A(a)$ is a $\pi(G_\mathcal U)$-invariant vector.
This implies that $b^A(a) = b^G(a) = 0$.
The cocycle condition then gives $b^A(ga) = b^A(a)$ for $g\in G_\mathcal U$ and $a\in A_\mathcal U$.
Hence $b^A$ factorizes through a cocycle $b'\colon Q_\mathcal U\to \H^A$.
This gives a controlled affine action of $Q_\mathcal U$ on $\H^A$.
Since $Q_\mathcal U$ has property (FH) it follows that there is an invariant vector $v_0\in \H^A$.
This vector satisfies $p^A(\sigma(g)v_0)) = v_0$ for all $g\in G_\mathcal U$.

Let $B\subseteq G_\mathcal U$ be a bounded set and $\epsilon > 0$.
Since $\sigma$ is well controlled, $\sigma(B)(v_0) \subseteq (p^A)^{-1}(v_0)$ is relatively compact, so we may select finitely many $v_1,\ldots,v_k\in (p^A)^{-1}(v_0)$ such that $\sigma(B)(v_0)$ is contained in a union of balls $\bigcup_{1\leq j\leq k} B(v_j,\tfrac\epsilon2)$.
We apply Lemma \ref{lem:L} to the representation $A_\mathcal U\curvearrowright \H^k$ and the vector $(v_1,\ldots,v_k)$.
We find $a_1,\ldots,a_n\in A$ and $t_1,\ldots,t_n\in [0,1]$ with sum 1 such that $\norm{\sum_{i=1}^nt_i\pi(a_i)v_j-v_0}<\tfrac\epsilon2$ for $1\leq j\leq k$.
Let $w = \sum_{i=1}^nt_i\sigma(a_i)(v_0)$.
Then for $g\in B$ we have $\norm{\sigma(g)(v_0) - v_j}<\tfrac\epsilon2$ for some $j$, and then
\begin{align*}
    \norm{\sigma(g)w-w} &= \norm{\sum_{i=1}^n(t_i\sigma(ga_i)(v_0)-t_i\sigma(a_i)(v_0))}\\
    &= \norm{\sum_{i=1}^n(t_i\sigma(a_i)\sigma(g)(v_0)-t_i\sigma(a_i)(v_0))}\\
    &= \norm{\sum_{i=1}^nt_i\pi(a_i)(\sigma(g)(v_0)-v_0)}\\
    &= \norm{\sum_{i=1}^nt_i\pi(a_i)(\sigma(g)(v_0)) - v_0}\\
    &\leq \norm{\sum_{i=1}^nt_i\pi(a_i)v_j - v_0} + \tfrac\epsilon2 \\
    &< \epsilon.
\end{align*}
So $w$ is an almost invariant vector for $B$.
This shows that any affine controlled representation of $G_\mathcal U$ has almost invariant vectors.

Now to finish the proof we use Proposition \ref{prop:ultralimitAffine}: if $G$ would not have property (FH), we find a controlled affine action of the ultraproduct $G_\mathcal U$ with positive displacement, and thus not with almost invariant vectors.
So $G$ must have property (FH) and hence property (T) since it has bounded geometry and its bornology is generated by a single set.
\end{proof}

\subsection{The group \texorpdfstring{$\Homeo^+_\mathbb Z(\mathbb R)$}{of lifts of homeomorphisms}}

Let $\Homeo^+_\mathbb Z(\mathbb R)$ be the group of lifts of orientation-preserving homeomorphisms of the unit circle, as described in \cite[\S3.4]{Ros21}.
It can be described explicitly as
\[\Homeo^+_\mathbb Z(\mathbb R) = \{f\colon\mathbb R\to\mathbb R\text{ homeomorphism}\mid f(x+n)=f(x)+n\text{ for all }x\in\mathbb R,n\in\mathbb Z\}.\]
Rosendal shows that the bornology $\OB$ can be described as follows: a subset $A\subseteq \Homeo^+_\mathbb Z(\mathbb R)$ is bounded exactly when the set $\{a(0) \mid a\in A\} \subseteq \mathbb R$ is bounded.
In particular it follows that $\Homeo^+_\mathbb Z(\mathbb R)$ is coarsely equivalent to $\mathbb Z$.
We can use Theorem \ref{thm:Serre} to prove that this group has property (T) as a bornological group.
\begin{corollary}\label{corInSection:HomeoZR}
 The group $\Homeo^+_{\mathbb Z}(\mathbb R)$ equipped with $\OB$ has coarse property (T). Hence, it also has coarse property FH (and in particular topological property FH).
 \end{corollary}
\begin{proof}
There is a short exact sequence of Polish groups:
\[1\to \Z\to\Homeo^+_{\Z}(\R)\to \Homeo^+(\R/\Z)\to 1,\]
where $\Z$ is central. 
Here $\Z$ has the discrete bornology and $\Homeo^+(\R/\Z)$ is bounded.
From the explicit description of the bornology on $\Homeo^+_\mathbb Z(\mathbb R)$ it is clear that $\Homeo^+(\R/\Z)$ is the image of a single bounded set.
Moreover, the bornology of $\Homeo^+_\mathbb Z(\mathbb R)$ is generated by the single bounded set $S = \{a\in \Homeo^*_\mathbb Z(\R) \mid 0\leq a(0)\leq 1\}$.

For $x\in \mathbb R$ denote by $\{x\}$ the fractional part of $x$.
Consider $a,b \in \Homeo^*_\mathbb Z(\R)$ given by $a(x) = x-\frac14$ and $b(x) = x+2\{x\}$ if $\{x\} \leq \frac14$ and $b(x) = x+\frac23-\frac23\{x\}$ if $\{x\}\geq \frac14$.
Then $[a,b](0) = \frac12$ and $[b^{-1},a^{-1}](\frac12) = 1$.
So $[b^{-1},a^{-1}][a,b]$ sends 0 to 1.
Then $\left([b^{-1},a^{-1}][a,b]\right)^n$ sends 0 to $n$.
So the collection $\{\left([b^{-1},a^{-1}][a,b]\right)^n\mid n\in\Z\}$ is coarsely dense in $\Homeo^+_\mathbb Z(\R)$, so in particular, the commutator subgroup is coarsely dense.
So the conditions of Theorem \ref{thm:Serre} apply, and $\Homeo^+(\R/\Z)$ has coarse property (T).
It also follows that this group has topological property FH, because this is a weaker statement.
\end{proof}
In the same chapter of his book, Rosendal considers two other similar examples of Polish groups:
\begin{itemize}
    \item  the group $\Aut_{\Z}(\mathbb Q)$ of order-preserving permutations
of $\mathbb Q$ commuting with integral translations, equipped with the permutation group topology.
\item the subgroup
$AC^*_
{\Z}(\R)$ of $\Homeo^+_{\Z}(\R)$ consisting of all $h$ so that $h$ and $h^{-1}$ are
absolutely continuous when restricted to $[0,1]$, equipped with a certain topology (see \cite{Her}).
\end{itemize}
For similar reasons, these two groups equipped with $\OB$ are quasi-isometric to $\Z$, and satisfy the assumptions of Theorem \ref{thm:Serre}, and so they have coarse property (T) as well.

\section{Bounded products}
In this section we prove Proposition \ref{prop-bornologies-bounded-product} and Theorems \ref{thm-bounded-products-topological-T} and \ref{thm-bounded-products-coarse-T}.

\begin{proposition}\label{propInSection-bornologies-bounded-product}
Let $\Gamma_n$ be a sequence of groups with finite sets of generators $S_n$.
Let $G = \prod_n\{G_n,S_n\}$ be their bounded product, with the bornology $\mathcal B$ and topology $\tau$ introduced in the Introduction.
Then we have $\mathcal B = \OB = \mathcal K$.
\end{proposition}
\begin{proof}
Let $A\in\mathcal B$.
We show that $A$ is relatively compact for $\tau$ by showing that every ultrafilter $\mathcal U$ on $\overline A$ has a limit point.
Since $A$ is bounded there is an integer $L$ such that for all $(a_n) \in A$ we have $a_n\in S_n^L$ for all $n$.
Then we can define the ultralimit $c_n = \lim_\mathcal Ua_n \in S_n^L$.
Then $(c_n) \in \overline A$.
We show that it is a limit of $\mathcal U$.
Let $(c_n)\cdot U_\rho$ be an open set containing $(c_n)$.
Then there is an $N$ such that $\rho(n) \geq 2L$ for all $n>N$.
Then $(c_n)\cdot U_\rho \cap \overline A \supseteq \{(a_n)\in \overline A \mid a_n = c_n \text{ for all }n\leq N\}$, and since the latter is in $\mathcal U$, so is the former, showing that $\mathcal U$ converges to $(c_n)$ and that $A\in\mathcal K$.

We always have $\mathcal K\subseteq \OB$ so it remains to show that $\OB\subseteq \mathcal B$.
We use Proposition \ref{prop-coarsely-bounded}(iii).
Suppose that $A\not\in\mathcal B$.
Then for any positive integer $i$ we can find $((a_i)_n) \in A$ and $n_i \in\mathbb N$ with $l((a_i)_{n_i}) > i2^i$.
Let $\rho\colon\mathbb N\to\mathbb Z_{\geq 1}$ be a proper function satisfying $\rho(n_i)\leq i$ for all positive integers $i$.
For $k\geq 1$ let $V_k = U_{2^k\rho}$.
Then $V_k^2\subseteq V_{k+1}$ and $\bigcup_kV_k = G$.
For all $i$ we have $((a_i)_n) \in A$, but $l((a_i)_{n_i}) > i2^i \geq 2^i\rho(n_i)$ so $((a_i)_n) \not\in V_i$.
So $A\not\subseteq V_i$, showing that $A\not\in\OB$.

\end{proof}

We recall the statement of Theorem \ref{thm-bounded-products-topological-T}:
\begin{theorem}\label{thmInSection-bounded-products-topological-T}
Let $\Gamma_n$ be a sequence of finite groups with symmetric generator sets $S_n$.
Suppose that $\sup_n|S_n|<\infty$.
Let $X$ be the coarse disjoint union of the Cayley graphs $\Cay(\Gamma_n,S_n)$.
Then the bounded product $G$ has topological property (T) if and only if $X$ is an expander sequence.
\end{theorem}
\begin{proof} Suppose $G$ has topological property (T).
    For each $n$, let $L^2_0\Gamma_n$ denote the Hilbert space of functions on $\Gamma_n$ with sum zero.
    The group $\Gamma_n$ acts by multiplication on $L^2_0\Gamma_n$.
    Precomposing with the projection $G\to\Gamma_n$ gives a representation of $G$ on $L^2_0\Gamma_n$.
    Taking the direct sum gives a representation of $G$ on $\bigoplus_nL^2_0\Gamma_n$.
    This is a continuous representation without any invariant vectors.
    By property (T), there are also no almost invariant vectors.
    This shows that $X$ is an expander sequence.
    
    Conversely, suppose that $X$ is an expander sequence.
    Let $M_n$ be the normalized adjacency matrix of $\Gamma_n$.
    There is a constant $h>0$ such that the spectrum $\sigma(M_n)$ is contained in $[h-1,1-h]\cup \{1\}$ for all $n$.
    Then the normalized adjacency matrix of the product group $\Gamma_1\times\cdots\times\Gamma_N$ is $M_1\tensor\cdots\tensor M_N$.
    Its spectrum is still contained in $[h-1,1-h]\cup\{1\}$.
    For any representation $\Gamma_1\times\cdots\times\Gamma_N \curvearrowright \H$, denote by $p_{\leq N}$ the projection on the invariant vectors.
    Then we get
    \begin{align*}
        h\norm{p_{\leq N}v-v} & \leq \norm{v-p_{\leq N}v} - \norm{(M_1\tensor\cdots\tensor M_N)(v-p_{\leq N}v)}\\
        &\leq \norm{(M_1\tensor\cdots\tensor M_N)v-v}\\
        &\leq \sup_{s\in S_1\times\cdots\times S_N}\norm{sv-v}.
    \end{align*}
    
    Now let $\pi\colon G\to U(\H)$ be a continuous representation with almost invariant vectors.
    Let $C = \prod_nS_n \subseteq G$.
    Let $k_0 = \lceil2h^{-1}\rceil$.
    Since the representation has almost invariant vectors, there is a unit vector $v_{k_0} \in \H$ such that $\norm{sv_{k_0}-v_{k_0}} \leq \frac1{k_0}$ for all $s\in C$.
    We will recursively construct unit vectors $v_k\in \H$ for $k\geq k_0$, with $\norm{sv_k-v_k} \leq \frac1k$ for $s\in C$, and $\norm{v_{k+1}-v_k}\leq \frac4{hk^2}$.
    Suppose $v_k$ exists.
    Since the representation is continuous, there is an integer $N$ such that for every $s\in \prod_{n>N}S_n$, we have $\norm{sv_k-v_k} < \frac1{k(k+1)}$.
    Let
    \[v_{k+1} = \frac{\left(\frac{k-2}k+\frac2kp_{\leq N}\right)v_k}{\norm{\left(\frac{k-2}k+\frac2kp_{\leq N}\right)v_k}}.\]
    We know that $\norm{p_{\leq N}v_k-v_k}\leq \frac1{hk}$, so $\norm{\left(\frac{k-2}k+\frac2kp_{\leq N}\right)v_k - v_k} \leq \frac2{hk^2}$.
    This also gives $\norm{\left(\frac{k-2}k+\frac2kp_{\leq N}\right)v_k} \geq 1-\frac2{hk^2}$, and $\norm{v_{k+1}-v_k}\leq \frac4{hk^2}$.
    
    Now let $s\in C$ and let $s_{>N}$ denote the projection of $s$ on $S_{N+1}\times S_{N+2}\times\cdots$.
    We get
    \begin{align*}
        \norm{sv_{k+1}-v_{k+1}} &= \frac{\frac{k-2}k(sv_k-v_k) + \frac2kp_{\leq N}(s_{>N}v_k-v_k)}{\norm{\left(\frac{k-2}k+\frac2kp_{\leq N}\right)v_k}}\\
        &\leq \frac{\frac{k-2}{k^2} + \frac2{k^2(k+1)}}{1-\frac2{hk^2}}\\
        &=\frac{k-2+\frac2{k+1}}{k^2-2h^{-1}}\\
        &= \frac{k^2-k}{k^2-2h^{-1}}\cdot \frac1{k+1}\\
        &\leq \frac1{k+1}.
    \end{align*}
    So $v_{k+1}$ satisfies the desired inequalities.
    
    Now the $v_k$ form a Cauchy sequence, so there is a limit $v$, and this will be an invariant vector.

\end{proof}
\begin{remark}
We could also consider the product topology on $G$ in this Theorem, which is weaker than $\tau$, because the representation constructed in the first part of the proof is also continuous for this topology.
\end{remark}

If we take all possible finite quotients (up to isomorphism) of $\Gamma$ (of which there are necessarily only countably many), we find a new criterion for property $(\tau)$.
\begin{corollary}
Let $\Gamma$ be a finitely generated residually finite group, generated by the symmetric set $S$.
Let $\Gamma_1,\Gamma_2,\ldots$ be all finite quotients of $\Gamma$ and let $G$ be their bounded product.
Then $\Gamma$ has property $(\tau)$ if and only if $G$ has topological property (T).
\end{corollary}

\begin{proof}
It is well known that $\Gamma$ has property $(\tau)$ if and only if the sequences of Cayley graphs $\Cay(\Gamma_n,\overline S)$ is an expander sequence \cite{Lub05}.
By Theorem \ref{thm-bounded-products-topological-T}, this is equivalent to $G$ having topological property (T).
\end{proof}

Now we turn to the coarse property (T) for the bounded product $G$.
For a sequence of finitely generated groups $(\Gamma_n,S_n)$ and an ultrafilter $\mathcal U$ on $\mathbb N$, we denote by $\prod_n\Gamma_n/\mathcal U$ the (bounded) ultraproduct, that is, the quotient of the set of sequences $G = \{(g_n)\in \Gamma_n \mid \sup_nl(g_n) < \infty\}$ by the equivalence relation $(g_n)\sim (h_n)$ if $\forall_\mathcal Un[g_n=h_n]$.
This is again a finitely generated group.
If the $S_n$ are bounded in size, we can identify them all with a fixed symmetric set $S$, and view $\Gamma_n$ as a quotient $\phi_n\colon F_S \twoheadrightarrow \Gamma_n$.
Then $\prod_n\Gamma_n/\mathcal U$ is equal to $F_S/\lim_\mathcal U\ker(\phi_n)$, where $\lim_\mathcal U\ker(\phi_n) = \{g\in F_S\mid\forall_\mathcal Un\colon g\in\ker(\phi_n)\}$.
For any quotient $\Gamma$ of $F_S$, we say the \emph{spectral gap} of $\Gamma$ is the largest $\gamma\in\mathbb R$ such that the spectrum of the operator $\Delta_S = \sum_{s\in S}1-s \in C^*_{\max}(\Gamma)$ is contained in $\{0\}\cup[\gamma,\infty)$.

Before proving Theorem \ref{thm-bounded-products-coarse-T} we need a lemma.
It is clear that any finite product of groups with property (T) has property (T) again.
The lemma quantifies this.
For a finitely generated group $\Gamma$, generated by a symmetric set $S$ containing the identity, let $M_S\in \mathbb C[\Gamma]$ denote the averaging operator $M_S=\frac1{|S|}\sum_{s\in S}s$.
Denote by $\sigma_{\max}(M_S)$ the spectrum of $M_S$ in the maximal $C^*$-algebra $C_{\max}^*(\Gamma)$.
Then $\Gamma$ has property (T) if and only if there is $\epsilon>0$ with $\sigma_{\max}(M_S)\subseteq[\epsilon-1,1-\epsilon]\cup\{1\}$.

\begin{lemma}\label{lem-product-(T)-quantitatively}
Let $\Gamma_1,\Gamma_2$ be finitely generated groups, with symmetric generating sets $S_1$ and $S_2$.
Suppose $\sigma_{\max}(M_{S_i})\subseteq [\epsilon-1,1-\epsilon]\cup\{1\}$ for $i\in\{1,2\}$.
Then the group $\Gamma_1\times\Gamma_2$ is generated by $S_1\times S_2$ and $\sigma_{\max}(M_{S_1\times S_2}) \subseteq [\epsilon-1,1-\epsilon]\cup\{1\}$.
\end{lemma}
\begin{proof}
Consider a representation of $\Gamma_1\times\Gamma_2$ on a Hilbert space $\H$ without constant vectors.
Let $v\in\H$.
Note that the subspace $\H^{\Gamma_1}$ of $\Gamma_1$-invariant vectors is a subrepresentation, and it contains no $\Gamma_2$-invariant vectors.
Let $w$ be the projection of $v$ onto $\H^{\Gamma_1}$.
Then $\norm{M_{S_1\times S_2}w} = \norm{M_{S_2}M_{S_1}w} = \norm{M_{S_2}w}\leq (1-\epsilon)\norm{w}$.
We also have $\norm{M_{S_1\times S_2}(v-w)} = \norm{M_{S_2}M_{S_1}(v-w)}\leq \norm{M_{S_1}(v-w)}\leq (1-\epsilon)\norm{v-w}$.
we get
\[\norm{M_{S_1\times S_2}v}^2 = \norm{M_{S_1\times S_2}w}^2 + \norm{M_{S_1\times S_2}(v-w)}^2 \leq (1-\epsilon)^2\norm{v}^2.\]
This shows that $\sigma_{\max}(M_{S_1\times S_2})\subseteq [\epsilon-1,1-\epsilon]\cup\{1\}$.
\end{proof}

We now recall the statement of Theorem \ref{thm-bounded-products-coarse-T}.
\begin{theorem}\label{thmInSection-bounded-products-coarse-T}
Let $\Gamma_n$ be a sequence of finite groups with symmetric generator sets $S_n$.
Suppose that $\sup_n|S_n|<\infty$.
Let $X$ be the coarse disjoint union of the Cayley graphs $\Cay(\Gamma_n,S_n)$.
The following are equivalent:
\begin{enumerate}
    \item The coarse space $X$ has geometric property (T).
    \item All ultraproducts $G_\mathcal U = \prod_n\Gamma_n/\mathcal U$ have property (T) with uniform spectral gap.
    \item All ultraproducts $G_\mathcal U$ have property (T).
    \item There is a finitely generated group $\Gamma$ with property (T) such that all $\Gamma_n$ are a quotient of $\Gamma$.
    \item The bounded product $G$ has coarse property (T).
    \item The bounded product $G$ has coarse property (T-).
\end{enumerate}
\end{theorem}
\begin{proof} We can assume without loss of generality that $|S_n|$ is constant. Then, 
we can view each $\Gamma_n$ as a quotient $\phi_n\colon F_S \twoheadrightarrow \Gamma_n$.
The ultralimit of these maps is $\phi_\mathcal U\colon F_S \twoheadrightarrow G_\mathcal U$.

Suppose (i) holds.
Let $\mathbb C_{\cs}[X]$ denote the Roe algebra of $X$, and let $\gamma > 0$ such that $\sigma(\pi(\Delta_S))\subseteq \{0\}\cup [\gamma,\infty)$ for all representations $\pi$.
Let $\mathcal U$ be an ultrafilter on $\mathbb N$ and consider a unitary representation $\pi\colon G_\mathcal U \to U(\H)$.
Suppose there is a unit vector $v\in \H$, and $\eta$ with $\pi(\Delta_S)v = \eta v$.
We will show that $\eta\in\{0\}\cup[\gamma,\infty)$, which implies that $G_\mathcal U$ has property (T) with spectral gap at least $\gamma$.

For any $n$, consider the Hilbert space $\ell^2(X_n,\H)$.
The group $F_S$ acts on this by $g\cdot \xi(x) = \pi(g)(\xi(g^{-1}x))$ for $\xi\in \ell^2(X_n,\H)$ and $x\in X_n$.
Now let $\H' = \prod_n\ell^2(X_n,\H)/\mathcal U$.
We define a representation $\pi'\colon\mathbb C_{\cs}[X] \to B(\H')$ as follows.
Any element of $\mathbb C_{\cs}[X]$ can (usually non-uniquely) be written as a finite sum $\sum_{g\in F_S}f_gg$ with $f_g\in\ell^\infty X$. In other words, $\mathbb C_{\cs}[X]$ is a quotient of the $*$-algebra $\ell^{\infty}(X)\rtimes F_S$.
We define $\pi'$ by 
\[\pi'\left(\sum_{g\in F_S}f_gg\right)((\xi_n)) = \left(\left(\sum_{g\in F_S}f_g\cdot (g\cdot \xi_n)\right)\right)\]
for $((\xi_n)) \in \prod_n\ell^2(X_n,\H)/\mathcal U$.
This is a well-defined $*$-representation of $\ell^{\infty}(X)\rtimes F_S$. To see that it induces a representation of $\mathbb C_{\cs}[X]$, suppose that $\sum_{g\in F_S}f_gg = 0$ in $\mathbb C_{\cs}[X]$.
For any $g,h \in F_S$ we know that if $\phi_n(g) = \phi_n(h)$ modulo $\mathcal U$, then $\pi(g) = \pi(h)$.
Since only finitely many $g$ occur in the sum, we know that there is $A\in \mathcal U$ such that for all $n\in A$ and all $g,h$ with $f_g\neq 0$ and $f_h\neq 0$, we have that $\phi_n(g) = \phi_n(h)$ implies $\pi(g) = \pi(h)$.
Now for $n\in A$ and $x\in X_n$ we have
\[\sum_{g\in F_S}f_g\cdot (g\cdot\xi_n) (x) = \sum_{y\in X_n}\sum_{g\mid gy = x}f_g(x)\pi(g)(\xi_n(y)).\]
Now if $g,h$ both satisfy $gy=x$ and $hy=x$ we have $\phi_n(g) = \phi_n(h)$, so $\pi(g) = \pi(h)$.
Moreover by the assumption that $\sum_{g\in F_S}f_gg = 0$ we know that $\sum_{g\mid gy=x}f_g(x) = 0$ for all $x,y\in X_n$, we see that the inner sum is zero for each $y$.
Therefore, the map $\pi'$ is well-defined representation of $\mathbb C_{\cs}[X]$.

Consider the unit vector $((\xi_n)) \in \H'$ given by $\xi_n(x) = |X_n|^{-\frac12}v$ for all $x\in X_n$.
Then $\pi'(\Delta_S)((\xi_n)) = \eta ((\xi_n))$.
This shows that $\eta \in \sigma(\pi'(\Delta_S)) \subseteq \{0\}\cup [\gamma,\infty)$, as we wanted to prove.

The implication $(ii)\to (iii)$ is trivial.
Now suppose $(iii)$.
First we prove that there is an integer $R$ such that the group $F_S/\langle\ker(\phi_n)\cap S^R\rangle$ has property (T), for all $n$.
If this is not the case, then for every $R$ there is an $n_R \in \mathbb N$ such that $F_S/\langle\ker(\phi_{n_R})\cap S^R\rangle$ does not have property (T).
Let $\mathcal U$ be a non-principal ultrafilter and consider the ultraproduct $H = \prod_R\Gamma_{n_R}/\mathcal U$.
This (discrete) group has property (T) by assumption.
By the main theorem of \cite{Oza16}, we can write
\[\Delta_S^2 - \gamma\Delta_S = \sum_{i=1}^k\xi_i^*\xi_i \in \mathbb R[H]\]
for some $\gamma > 0$ and $\xi_i \in \mathbb R[H]$.
Let $\tilde\xi_i$ be lifts of $\xi_i$ in $\mathbb R[F_S]$.
For large enough $R$, the supports of the $\tilde\xi_i$ will be contained in $S^{R/4}$.
Now for most $R$ (in the sense of the ultrafilter) we have
\[\Delta_S^2-\gamma\Delta_S-\sum_{i=1}^k\xi_i^*\xi_i \in \ker(\mathbb R[F_S] \twoheadrightarrow \mathbb R[\Gamma_{n_R}]) \cap \mathbb R[S^{R/2}] \subseteq \ker(\mathbb R[F_S] \twoheadrightarrow \mathbb R[F_S/\langle\ker(\phi_{n_R})\cap S^R\rangle]).\]
So we see that for these $R$ the spectral gap of $F_S/\langle \ker(\phi_{n_R})\cap S^R\rangle$ is at least $\gamma$, giving a contradiction.

So let $R$ be such that $F_S/\langle\ker(\phi_n)\cap S^R\rangle$ has property (T) for all $n$.
Let $\Lambda_1,\ldots,\Lambda_m$ be all the groups that can arise in this way.
There are only finitely many since there are only finitely many subsets of $S^R$.
All these groups have property (T).
Then their product $\Gamma$ has property (T), and all $\Gamma_n$ are a quotient of $\Gamma$, showing the implication $(iii)\to(iv)$.

Now suppose $(iv)$ is true, so all $\Gamma_n$ are quotients of a single group $\Gamma$ with property (T).
We can also assume that $\Gamma$ is generated by $S$ and $S_n$ is the image of $S$.
We will show that $G$ has coarse property FH, which is equivalent to coarse property (T) since the bornology of $G$ is generated by the bounded set $\prod_nS_n$.

We can assume that $S$ is symmetric and contains the identity.
Then there is $\epsilon > 0$ such that $\sigma_{\max}(M_S) \subseteq [\epsilon-1,1-\epsilon]\cup\{1\}$.
Let $\alpha\colon G\curvearrowright \H$ be a controlled affine action.
Let $v\in\H$ and let $C>0$ such that $\norm{\alpha(s)v-v}\leq C$ for all $s\in \prod_nS_n$.
Let $g\in G$.
We will show that $\norm{\alpha(g)v-v}\leq 2C/\epsilon$, showing that $\alpha$ has bounded orbits.

Since $g$ is in the bounded product, there is an integer $k$ such that $g\in (\prod_nS_n)^k$.
Each coordinate of $g$ is then an image of an element in $S^k$.
Since $S^k$ is finite, this means there is an integer $m$ and a homomorphism $f\colon \Gamma^m\to G$ and an element $x\in \Gamma^m$ with $f(x) = g$.
By Lemma \ref{lem-product-(T)-quantitatively}, we have $\sigma_{\max}(M_{S^m}) \subseteq [\epsilon-1,1-\epsilon]\cup\{1\}$.

Let $v_0$ denote the orthogonal projection of $v$ on the subspace of $\H$ consisting of $\Gamma^m$-fixed vectors.
Define $\pi\colon \Gamma^m \curvearrowright \H$ by $\pi(y)w = \alpha(f(y))(w+v_0)-v_0$.
This defines a representation of $\pi$ on $\H$.
Since $v-v_0$ is perpendicular to the constant vectors of this representation, we have $\norm{\pi(M_{S^n})(v-v_0)} \leq (1-\epsilon)\norm{v-v_0}$.
We also have
\[\norm{\pi(M_{S^n})(v-v_0)} \geq \norm{v-v_0} - \norm{\frac1{|S^m|}\sum_{s\in S^m}\alpha(f(s))v-v}\geq \norm{v-v_0}-C.\]
This gives $\norm{v-v_0} \leq C/\epsilon$.
Finally we have
\[\norm{\alpha(g)v-v} = \norm{\alpha(f(x))v-v} = \norm{\pi(x)(v-v_0)+v_0-v}\leq 2\norm{v_0-v} \leq 2C/\epsilon.\]

This shows that $\alpha$ has bounded orbits, thus $G$ has coarse property FH, and coarse property (T) as well, showing $(iv)\to(v)$.

Suppose $G$ has coarse property (T).
We will prove that $X$ has geometric property (T).
Let $(A,\epsilon)$ be a Kazhdan pair for $G$.
Let $N = \sup_{(g_n)\in A}\sup_nl(g_n)$ which is finite since $A$ is bounded.
Let $\pi\colon\mathbb C_{\cs}[X] \to B(\H)$ be a representation.
Let $v\in \H_c^\perp$ be a unit vector, perpendicular to all constant vectors.
Note that $G$ acts on $L^2X$, and this gives a natural embedding $G\to \mathbb C_{\cs}[X]$.
Thus $\pi$ restricts to a representation of $G$.
The constant vectors form a subrepresentation, so $\H_c^\perp$ is also a representation of $G$.
It does not have any constant vectors.
So there is $a\in A$ with $\norm{av-v}>\epsilon$.
Now $a$ can be viewed as an element in $\mathbb C_{cs}[X]$ with propagation at most $N$.
So it follows that $X$ has geometric property (T), showing $(v)\to(i)$.

The implication $(v)\to(vi)$ is trivial.
We will finish by showing $(vi)\to(iii)$.
Suppose $G$ has property (T-), or equivalently FH-, and let $\alpha\colon G_\mathcal U\curvearrowright \H$ be an affine action.
Then the composition $G\to G_\mathcal U \to U(\H)$ gives a well-controlled action, because the image of any bounded subset of $G$ is finite.
This shows that a fixed point exists, so $G_\mathcal U$ has property FH and property (T).
\end{proof}

\section{Further remarks and questions}\label{sec:questions}

\noindent{\bf Cohomological interpretation.}
By the Mazur-Ulam theorem, an affine isometric action of $G$ on a Banach space $E$ corresponds a pair $(\pi,b)$, where $\pi$ is a norm-preserving representation
$\pi$ of $G$ on $E$, and  $b\in Z^1(G,\pi)$  is $1$-cocycle, i.e.\ a map $b:G\to E$ satisfying the cocycle relation $b(gh)=\pi(g)b(h)+b(g)$ for all $g,h\in G$. The subspace $B^1(G,\pi)$ of $Z^1(G,\pi)$  consists of coboundaries, i.e.\ cocycles of the form $b(g)=v-\pi(g)v$ where $v\in E$ (this corresponds to the affine action fixing the point $v$).
The first cohomology groups with values in $\pi$ is defined as $H^1(G,\pi)=Z^1(G,\pi)/B^1(G,\pi)$. 

Assume that $G$ is a bornological group. Note that an affine action given as $\sigma(g)v=\pi(g)v+b(g)$. Hence $\sigma$ is controlled if and only if $b$ is ``controlled'': i.e.\ bounded on bounded sets. The subspaces of controlled cocycles will be denoted by $Z^1_{\born}(G,\pi)$. Note that $B^1(G,\pi)\subset Z^1_{\born}(G,\pi)$. The {\it controlled first cohomology} will be defined as
$H^1_{\born}(G,\pi)=Z^1_{\born}(G,\pi)/B^1(G,\pi)$.

If $G$ is a topological group, and $\pi$ is a norm-preserving representation (not necessarily continuous), we let $Z^1_{\cont}(G,\pi)$ be the subspace of continuous $1$-cocycles, and let $H^1_{\cont}(G,\pi)=Z^1_{\cont}(G,\pi)/B^1(G,\pi)$.
If $G$ is locally compact, equipped with the bornology $\mathcal K$, then the inclusion 
 $Z^1_{\cont}(G,\pi)\to Z^1_{\born}(G,\pi)$ gives rise to an injective linear map $H^1_{\cont}(G,\pi)\to H^1_{\born}(G,\pi)$.

\begin{theorem}\label{thm:cohomo}
Let $G$ be a locally compact group. Let $\pi$ be a norm-preserving (not necessarily continuous) representation on a uniformly convex Banach space $E$. Then every controlled $1$-cocycle is cohomologous to a continuous one. In other words, the map
$H^1_{\cont}(G,\pi)\to H^1_{\born}(G,\pi)$ is an isomorphism.
\end{theorem}
\begin{proof}
Let $b\in Z^1_{\born}(G,\pi)$, and let $\sigma$ be the corresponding (controlled) affine action. Proposition \ref{prop:continuousSubaction} implies that there is $v\in E$ such that $g\mapsto\sigma(g)\cdot v$ is continuous. Let $c$ be the $1$-coboundary defined by $c(g)=\pi(g)v-v$. Then the cocycle $b'=b+c$ satisfies 
$b'(g)=b(g)+\pi(g)v-v=\sigma(g)\cdot v-v$. It is therefore continuous. 
\end{proof}

Note that the previous definitions extend to higher degree cohomology groups. In particular, we have a injective map $H^k_{\cont}(G,\pi)\to H^k_{\born}(G,\pi)$.

\begin{que}
Under the assumptions of Theorem \ref{thm:cohomo},
is the map $H^k_{\cont}(G,\pi)\to H^k_{\born}(G,\pi)$ surjective (hence bijective) for all $k\geq 1$?
\end{que}

\

\noindent{\bf Shopping centres for Banach spaces.}
The proofs of  Theorems \ref{thm-(T)-and-(T-)-equivalent} and \ref{thm:Serre} rely on the notion of shopping centre (see \S \ref{sec:shoppingcentres}) which attach to every bounded subset of a Hilbert space a relatively compact subset which satisfies a crucial property stated in Lemma \ref{lem-mall-relatively-compact}. We have not been able to define an analogous notion for uniformly convex Banach spaces. 

\begin{que}
Is there an analogous version of Lemma \ref{lem-mall-relatively-compact} for uniformly convex Banach spaces? 
\end{que}

\begin{que}
Are Theorems \ref{thm:Serre} and \ref{thm-(T)-and-(T-)-equivalent} valid for uniformly convex Banach spaces? 
\end{que}

\

\noindent{\bf Bounded products.}
The notion of bounded product makes sense for a family of groups equipped with a left-invariant pseudo-metric. It would be interesting to investigate it more systematically (see also \cite{Win21}). In particular:
\begin{que}
Can Theorems \ref{thm-bounded-products-topological-T} and \ref{thm-bounded-products-coarse-T} have analogues for sequences of compact groups equipped with suitable left-invariant metrics? 
\end{que}

\

\noindent{\bf Stability of coarse FH.}
Among locally compact groups, property FH is invariant under taking closed cocompact unimodular subgroups. In particular this holds if the subgroup is closed cocompact and normal.
    The unimodulaityr assumption is crucial as for example, $\textnormal{SL}_3(\mathbb R)$ has property FH, but the closed cocompact subgroup of upper-triangular matrices does not (as it is solvable). 
    The situation is even worse among Polish groups: indeed we have seen that $\Homeo^+_\mathbb Z(\mathbb R)$ has (even coarse) property FH but contains $\Z$ as a coarsely dense discrete subgroup. Note that for locally compact groups, the unimodularity assumption of the subgroup $H<G$ is here to ensure that $G/H$ have a $G$-invariant probability measure.  
     \begin{que}
Is there a reasonably natural restriction on pairs $(G,H)$, where $G$ is a bornological group and $H$ is a subgroup which ensures that property FH for $G$ implies property FH for H?

\end{que}

\vfill
\centering
\small

Romain Tessera\\
\textsc{University of Paris Cité, Sorbonne Université,}\\
 Institut de Mathématiques de Jussieu-Paris Rive Gauche, F-75013 Paris, France\\
romatessera@gmail.com\\
\vspace{2cm}
Jeroen Winkel\\
\textsc{Westf\"alische Wilhelms-Universit\"at M\"unster}\\
Einsteinstrasse 62, 48149 M\"unster\\
Winkeljeroen@gmail.com
\end{document}